\documentclass[11pt]{amsart}

\usepackage[a4paper,margin=1in]{geometry}
\usepackage{amsmath,amssymb,amsthm,mathtools}
\usepackage{microtype}
\usepackage{booktabs}
\usepackage{graphicx}
\usepackage{float}

\usepackage[numbers,sort&compress]{natbib}
\usepackage[colorlinks=true,linkcolor=blue,citecolor=blue,urlcolor=blue]{hyperref}

\numberwithin{equation}{section}

\newtheorem{theorem}{Theorem}[section]
\newtheorem{proposition}[theorem]{Proposition}
\newtheorem{lemma}[theorem]{Lemma}

\newtheorem{remark}[theorem]{Remark}

\newcommand{\R}{\mathbb{R}}
\newcommand{\Pp}{\Pi_+}
\newcommand{\dd}{\,\mathrm{d}}
\newcommand{\ind}{\mathbf{1}}
\newcommand{\Ecal}{\mathcal{E}}

\title[Linear Growth of the Vorticity Maximum]
{ Linear Growth of the Vorticity Maximum for Axisymmetric Euler Flows Without Swirl}

\author{Daomin Cao}
\address{State Key Laboratory of Mathematical Sciences, Academy of
Mathematics and Systems Science, Chinese Academy of Sciences, Beijing
100190, P.R. China}
\email{dmcao@amt.ac.cn}

\author{Junhong Fan}
\address{Institute of Applied Mathematics, Academy of Mathematics and
Systems Science, Chinese Academy of Sciences, Beijing 100190, P.R. China;
and University of Chinese Academy of Sciences, Beijing 100049,
P.R. China}
\email{fanjunhong@amss.ac.cn}

\author{Guolin Qin}
\address{State Key Laboratory of Mathematical Sciences, Academy of
Mathematics and Systems Science, Chinese Academy of Sciences, Beijing
100190, P.R. China;
and University of Chinese Academy of Sciences, Beijing 100049,
P.R. China}
\email{qinguolin18@mails.ucas.ac.cn}

\date{}

\hypersetup{
  pdftitle={Linear Growth of the Vorticity Maximum for Axisymmetric Euler Flows Without Swirl},
  pdfauthor={Daomin Cao, Junhong Fan, Guolin Qin},
  pdfsubject={Axisymmetric Euler equations without swirl},
  pdfkeywords={Euler equations, axisymmetric flow, radial moment, vortex stretching, vortex-ring collision, vortex patch, linear vorticity growth}
}

\begin{document}
\raggedbottom

\begin{abstract}
We study long-time vortex stretching for three-dimensional
axisymmetric Euler flows without swirl in the anti-parallel class
associated with the head-on collision of two coaxial vortex rings.
This geometry motivated Childress's \(t^{4/3}\) conjecture for the
vorticity maximum in the full axisymmetric no-swirl class
[S.~Childress, \emph{Physica D} \textbf{237} (2008), 1921--1925].
For unit-strength relative-vorticity patches in this class, we prove
that the outer radius, which is exactly the vorticity maximum, reaches
the linear scale on an arbitrarily large fixed proportion of every
sufficiently large dyadic interval:
for every \(0<\eta<1\), there exist \(c_\eta>0\) and \(T_\eta>1\) such
that
\[
 \left|
 \left\{t\in[T,2T]:
 \mathcal R_\omega(t)
 =\|\boldsymbol\Omega(t)\|_{L^\infty(\mathbb R^3)}
 \ge c_\eta t
 \right\}
 \right|\ge(1-\eta)T
 \qquad(T\ge T_\eta).
\]
The same estimate for \(\mathcal R_\omega(t)\) holds for all data
considered below.
For every nontrivial compactly supported initial datum in this class
that is odd in \(z\) and non-positive for \(z>0\), we also prove
\[
 \lim_{t\to\infty}
 {[\log(2+t)]^{5/2}(1+t)^{-3/2}\iint_{\Pi_+}r^2[-\omega(r,z,t)]\dd r\dd z}=+\infty.
\]
To the best of our knowledge, this is the first radial-moment lower bound with exponent greater than
one.  For unit-strength patches, the same moment bound also yields the
full-time estimate
\[
 \lim_{t\to\infty}
 \frac{\|\boldsymbol\Omega(t)\|_{L^\infty(\mathbb R^3)}
 [\log(2+t)]^{5/4}}{(1+t)^{3/4}}
 =+\infty.
\]
For general data, we further obtain a quantitative Eulerian form of
simultaneous radial escape and collision.  The proof uses two monotone
mixed moments, a compactly supported multiplier, and an exterior
\(L^2\) estimate for the velocity generated by interior vorticity.
\end{abstract}

\maketitle
\enlargethispage{4pt}

\section{Introduction and main results}\label{sec:intro}

Quantifying the growth of vorticity in a globally regular
three-dimensional Euler flow is a basic problem in fluid dynamics.  The
velocity and pressure solve the incompressible Euler equations
\[
 \partial_t u+(u\cdot\nabla)u=-\nabla p,\qquad \nabla\cdot u=0,
\]
and the vorticity \(\boldsymbol\Omega=\nabla\times u\) satisfies
\begin{equation}\label{eq:intro-vorticity}
 \partial_t\boldsymbol\Omega+(u\cdot\nabla)\boldsymbol\Omega
 =(\boldsymbol\Omega\cdot\nabla)u,\qquad
 u=\nabla\times(-\Delta)^{-1}\boldsymbol\Omega.
\end{equation}
The stretching term on the right-hand side of \eqref{eq:intro-vorticity}
vanishes in two dimensions and is the essential three-dimensional
mechanism.  The Beale--Kato--Majda criterion makes the vorticity maximum
central to the regularity problem \cite{BKM}; even within classes of
global solutions, however, its possible long-time growth is not yet
fully understood.  We refer to
\cite{MajdaBertozzi,KiselevReview,DrivasElgindi} for broader
accounts of vorticity growth and singularity formation.

Axisymmetric flows without swirl provide a particularly transparent
setting for this question. We will use the cylindrical coordinate
frame $\{\mathbf{e_r},\mathbf{e_z},\mathbf{e_\theta}\}$.  If
\(\boldsymbol\Omega=\omega(r,z,t)\mathbf{e_\theta}\), then the relative
vorticity \(\omega/r\) is transported by the meridional flow.  Thus
\(\omega/r\) does not amplify along trajectories, whereas
\(\omega=r(\omega/r)\) grows when vorticity-bearing particles move
outward.  Radial transport is therefore not merely a geometric effect:
it is precisely the vortex-stretching mechanism in this class.  Global
existence and regularity under several smoothness, decay, and
integrability assumptions were developed in
\cite{UkhovskiiIudovich,Majda,SaintRaymond,ShirotaYanagisawa,Serfati}
and refined in \cite{DanchinPerfect,Danchin,AbidiHmidiKeraani};
axisymmetric vortex patches were treated in
\cite{GamblinSaintRaymond}.

At lower H\"older regularity, the same axisymmetric no-swirl structure
can instead develop finite-time singularities.  Elgindi constructed
such solutions for \(C^{1,\alpha}\) velocity fields with sufficiently
small \(\alpha>0\) \cite{Elgindi}; Elgindi--Ghoul--Masmoudi established
stability and obtained finite-energy solutions in the same regularity
regime \cite{ElgindiGhoulMasmoudi}.  A different finite-energy
construction, again for sufficiently small \(\alpha>0\), with velocity in
\(C^\infty(\R^3\setminus\{0\})\cap C^{1,\alpha}\cap L^2(\R^3)\)
was given in \cite{CordobaMartinezZoroaZheng}.  Independent recent
preprints extend finite-time blow-up constructions to every
\(0<\alpha<1/3\)
\cite{Shkoller,ChenBlowup,ShaoWeiZhangZhangBlowup}.
The present paper concerns the complementary question of infinite-time
growth in the globally regular class specified below.

Childress \cite{ChildressReports,Childress} proposed \(t^{4/3}\) as the
natural candidate scale for
\(\|\boldsymbol\Omega(t)\|_{L^\infty(\R^3)}\).  The model behind the
conjecture is a pair of coaxial vortex rings of opposite circulation,
placed symmetrically above and below \(z=0\).  In a meridional section
they appear as a counter-rotating pair.  Their mutual induction has two
coupled effects: the axial motion brings the rings toward the symmetry
plane, while the active vorticity is driven radially outward.  Since
\(\omega/r\) is transported, the increase of \(r\) lengthens the
circular vortex lines and amplifies \(\omega\) by the same factor.
For a material vortex tube, incompressibility simultaneously forces
its transverse area to contract.  The conjecture concerns the full axisymmetric no-swirl
class; the odd, one-signed class considered here is the anti-parallel
collision geometry expected to produce large growth.

In the formal picture of Childress and
Childress--Gilbert--Valiant \cite{ChildressGilbertValiant}, the active
vorticity organizes into a thin expanding dipole head whose rescaled
cross-section is close to a Sadovskii dipole.  The head approximately
retains its kinetic energy while circulation may be shed into a
trailing filament.  This loss of circulation is the erosion mechanism.
Childress--Gilbert--Valiant developed the corresponding asymptotic
model and verified its predictions numerically.  The resulting picture
is not simply one of radial spreading: outward escape, axial approach,
transverse compression, vorticity amplification, and possible erosion
are parts of the same dynamics.

Head-on collisions have also been studied experimentally and
computationally, mainly for viscous flows, in
\cite{KambeMinota,StanawayShariffHussain,LimNickels,ChuEtAl,GuanEtAl};
see also the review \cite{ShariffLeonard}.  In viscous experiments and
simulations, flattening of the rings may be followed by reconnection,
annihilation, or acoustic emission.  Such effects are absent from the
inviscid model considered here, but the common geometric event is the
approach of opposite rings to a symmetry plane together with rapid
radial expansion.  A rigorous construction of an Euler solution
attaining the predicted \(t^{4/3}\) growth remains open.

On the upper-bound side, the \(t^{4/3}\) scale has now been proved in
several important but distinct classes of data.  Lim--Jeong
\cite{LimJeong} treated compactly supported vorticity.
Shao--Wei--Zhang \cite{ShaoWeiZhang} and Egamberganov--Yao
\cite{EY} obtained the same rate without compact-support assumption, under
different regularity or weighted integrability hypotheses.  These results establish
the conjectured \(t^{4/3}\) upper bound within their respective classes;
they do not show that this rate is attained.

The lower-bound problem is substantially harder.  Transport of
\(\omega/r\) explains how radial motion creates stretching, but it does
not force that motion to occur.  A convenient measure of radial escape
is
\begin{equation}\label{eq:intro-P}
 P(t):=\iint_{\Pp}r^2[-\omega(r,z,t)]\,\dd r\dd z,
\end{equation}
where the sign convention is chosen so that \(-\omega\geq0\) in the
upper half-plane
\(\Pp=\{(r,z):r>0,\ z>0\}\).  Choi--Jeong proved monotonicity formulas for
the radial and vertical moments \cite[Theorem~1.1]{CJ} and obtained
\(P(t)\gtrsim_\varepsilon t^{2/15-\varepsilon}\)
\cite[Theorem~1.3]{CJ}.  Gustafson--Miller--Tsai combined this
framework with positive-kernel energy estimates and improved the bound
to \(P(t)\gtrsim_\varepsilon t^{3/4-\varepsilon}\) for the smooth
solutions in their class
\cite[Theorem~1.4 and Sections~3.1--3.3]{GMT}.  Most recently,
Egamberganov--Yao derived the exact axis identity
\begin{equation}\label{eq:intro-EY-axis}
 P'(t)=\int_0^\infty r\,u^r(r,0,t)^2\,\dd r
\end{equation}
and proved
\[
 P(t)\gtrsim\frac{t}{\log(2+t)}
\]
under their sign, symmetry, integrability, and moment assumptions
\cite[Theorem~1.2 and Section~4]{EY}.  They also obtained
\(t^{1/4}\) limsup growth of every vorticity \(L^p\)-norm,
\(1\leq p\leq\infty\), and quantitative time-integrated escape from
each fixed radial cylinder \cite[Theorems~1.4 and~1.5]{EY}.  On the
upper-bound side, they established
\[
 P(t)\lesssim(1+t)^2,
\]
the rate predicted to be optimal in Childress's scenario
\cite[Theorem~1.1]{EY}.  The present three-halves lower bound therefore
advances the attainability problem for this conjectured quadratic
radial-moment scale.  For unit-strength patches, we additionally reach
the genuinely linear scale for the vorticity maximum on large sets of
times.

The proofs of the three preceding lower bounds, those of
Choi--Jeong, Gustafson--Miller--Tsai, and Egamberganov--Yao, all use
the decreasing vertical moment
\[
 Z(t):=\iint_{\Pp}z[-\omega(r,z,t)]\,\dd r\dd z
\]
as an essential a priori input; see \cite[Lemma~3.3]{CJ},
\cite[Section~3.1]{GMT}, and \cite[Sections~4--5]{EY}.  Its
dissipation controls \(u^r\) in the bulk and \(u^z\) on the axis, but
the weight \(z\) does not couple axial position to radial escape.  Our
proof does not use \(Z\).  Instead, we use two mixed radial--axial
moments defined in \eqref{MPhi} and \eqref{eq:critical-Psi}, respectively.  The first gives
\[
 \sup_{t\geq0}\iint_{\Pp}rz[-\omega(r,z,t)]\,\dd r\dd z
 \leq C(\omega_0),
\]
and therefore, for every \(R,H>0\),
\[
 \iint_{\{r\geq R,\ z\geq H\}}[-\omega(r,z,t)]\,\dd r\dd z
 \leq \frac{C(\omega_0)}{RH}.
\]
The bound obtained from \(Z\) alone has no factor \(R^{-1}\).
Thus the axial control becomes stronger as the vorticity moves to large
cylindrical radius.  Combined with the energy estimates, this locates
the energy-generating escaping vorticity in a shrinking neighborhood of
\(z=0\).  The dissipation of the second moment controls both independent
components of the meridional velocity with the scale-dependent weights
needed to localize the conserved kinetic energy.  These two estimates are
the main new ingredients behind the three-halves lower bound for \(P(t)\),
the linear support estimate, and the corresponding patch conclusion.

The progression of lower bounds is summarized in Table~\ref{tab:lower-bounds}.
The notation \(t^{\alpha-}\) means every power strictly below \(t^\alpha\). For simplicity, we have ignored some more precise logarithmic terms.

\begin{table}[H]
\centering
\caption{Lower bounds for the radial moment and the vorticity maximum}
\label{tab:lower-bounds}
\small
\renewcommand{\arraystretch}{1.38}
\begin{tabular}{@{}p{.27\textwidth}p{.27\textwidth}p{.40\textwidth}@{}}
\toprule
Reference
& radial moment \(P(t)\)
& vorticity maximum
  \(\displaystyle\|\boldsymbol\Omega(t)\|_{L^\infty(\mathbb R^3)}\)\\
\midrule
Choi--Jeong \cite{CJ}
& \(t^{2/15-}\)
& \(t^{1/15-}\)\\

\raggedright Gustafson--Miller--Tsai \cite{GMT}
& \(t^{3/4-}\)
& \(t^{3/8-}\)\\

Egamberganov--Yao \cite{EY}
& \(t^{-}\)
& \(t^{1/2-}\)\\

Present work
& \(\displaystyle t^{3/2-}\)
& \(\displaystyle t^{3/4-}\)\\

Present work
& ---
& \(\displaystyle t\) on a proportion \(1-\eta\)
  for every fixed \(0<\eta<1\)\\
\bottomrule
\end{tabular}
\end{table}

All rows except the last record full-time pointwise estimates.  The
patch entries in the first four rows follow from the corresponding
moment estimates; see
\cite[Corollary~1.4 and Remark~1.5]{CJ},
\cite[Remark~1.6]{GMT}, and
\cite[Corollary~1.3]{EY}.
The main theorem of Gustafson--Miller--Tsai is formulated for sufficiently
smooth solutions, while their Remark~1.6 records the moment-to-patch
implication displayed in the table.

The principal conclusion for the vorticity itself is genuinely linear.
For every unit-strength relative-vorticity patch considered below and
every \(0<\eta<1\), there exist \(c_\eta>0\) and \(T_\eta>1\) such that
\[
 \left|
 \left\{t\in[T,2T]:
 \mathcal R_\omega(t)
 =\|\boldsymbol\Omega(t)\|_{L^\infty(\mathbb R^3)}
 \ge c_\eta t
 \right\}
 \right|\ge(1-\eta)T
 \qquad(T\ge T_\eta).
\]
Thus, after any proportion strictly below one has been fixed, the
vorticity maximum is bounded below by a positive multiple of \(t\) on
at least that proportion of every sufficiently large dyadic interval.
The estimate for \(\mathcal R_\omega(t)\) does not use the patch
structure and holds for every solution in Theorem~\ref{thm:main}.

For general data, we also obtain a radial-moment exponent greater than
one.  For all data in the class specified below,
\begin{equation}\label{eq:intro-main-informal}
 \frac{P(t)[\log(2+t)]^{5/2}}{(1+t)^{3/2}}
 \longrightarrow\infty.
\end{equation}
To the best of our knowledge, this is the first lower bound for \(P(t)\) with a
polynomial exponent larger than one.  We stress that \(P(t)\)
is a radial moment and not the vorticity maximum in Childress's
conjecture; the two quantities have different scaling.  Nevertheless,
\eqref{eq:intro-main-informal} gives a robust quantitative form of
radial spreading without assuming that a coherent Childress head or a
Sadovskii limiting profile has formed.

The proof yields a substantially stronger conclusion than the
support statement alone.  Theorem~\ref{thm:escape-collision} shows
that vorticity in a radially escaping region and a shrinking
neighborhood of the symmetry plane generates a velocity field carrying
an arbitrarily large fraction of the conserved kinetic energy.  The escape therefore
cannot be caused solely by an energetically negligible filament.  By
odd symmetry, there is an oppositely signed reflected part below the
plane.  The theorem thus captures simultaneous radial escape and
asymptotic approach to the collision plane in a precise energetic,
Eulerian sense.  It does not track a single coherent head and is
compatible with simultaneous erosion.  The mechanism requires no
patch assumption, boundary regularity, or stability of a prescribed
profile.

A central feature of the argument is that it is not restricted to
patches.  For general compactly supported data in our class,
Theorem~\ref{thm:general-Lp} proves growth of every vorticity
\(L^p\)-norm, uniformly over \(1\leq p\leq\infty\).  The full-time
estimate has the scale
\((1+t)^{1/4}[\log(2+t)]^{-25/12}\) with a divergent prefactor.  In
particular, it upgrades the \(t^{1/4}\) limsup phenomenon of
Egamberganov--Yao to a pointwise-in-time statement, up to logarithms.
On a density-one subset of each dyadic interval, the larger scale
\(T^{1/2}(\log T)^{-3/2}\) holds simultaneously for all \(p\).

For vortex patches, radial escape translates more directly into
vorticity amplification.  For the unit-strength patches considered
below, \(-\omega=r\) on the upper patch, so the vorticity maximum is
exactly its outer radius.  The linear support estimate gives the
principal vorticity conclusion stated above.  The radial-moment bound
also gives the full-time pointwise estimate
\[
 \frac{\|\boldsymbol\Omega(t)\|_{L^\infty(\R^3)}
 [\log(2+t)]^{5/4}}{(1+t)^{3/4}}\longrightarrow\infty,
\]
improving the preceding full-time polynomial exponent \(1/2\).  The
linear estimate is not asserted at every large time, and Childress's
\(t^{4/3}\) conjecture remains open.

For context, a distinct but closely related body of work concerns the
construction and dynamics of vortex rings.  It includes clustered
travelling rings \cite{AoLiuWei}, smooth leapfrogging rings
\cite{DavilaDelPinoMussoWei}, scaling limits or effective dynamics for
concentrated rings \cite{ButtaCavallaroMarchioro,DonatiEtAl}, and the
recent global-in-time concentration and propagation theorem for a thin
single vortex ring \cite{GuoJeongZhao}; see also Dyson's classical
calculation \cite{Dyson}.  Singular initial rings and their viscous
evolution were studied in \cite{FengSverak,GallaySverak}.  The
two-dimensional Sadovskii patch suggested by the rescaled Childress
head was constructed independently in
\cite{CJSadovskii,HuangTong}.  A scaling-invariant variational family
ranging from the Sadovskii patch to axis-touching vortices of
arbitrarily high finite regularity was subsequently obtained in
\cite{AbeChoiJeongSimWoo}.  Related filamentation and
vorticity-gradient growth appear in \cite{GuoJeongZhao,CJHill,Do}.  At lower
regularity, global weak solutions were constructed in
\cite{JiuLiuNiu}; renormalization and energy conservation were studied
in \cite{NobiliSeis}.  Recent norm-inflation results at the Lorentz-space
endpoint of the assumptions in \cite{Danchin} are given in
\cite{BangCheskidov}.

\subsection{Setting and main results}
Let
\[
 \Pi=\{(r,z):r>0,\ z\in\R\},\qquad
 \Pp=\{(r,z):r>0,\ z>0\}.
\]
We write
\[
 u=u^r(r,z,t)\mathbf{e_r}+u^z(r,z,t)\mathbf{e_z},\qquad
 \boldsymbol\Omega=\nabla\times u=\omega(r,z,t)\mathbf{e_\theta},\qquad
 \omega=\partial_z u^r-\partial_r u^z.
\]
The divergence-free condition gives 
$\partial_r(r u^r)+\partial_z(r u^z)=0$.
We assume that the initial data \(\omega_0\not\equiv0\) is compactly supported and
satisfies
\begin{equation}\label{eq:symmetry-sign}
 \omega_0(r,-z)=-\omega_0(r,z),\qquad
 \omega_0(r,z)\leq0\quad(z>0),
\end{equation}
together with
\begin{equation}\label{eq:admissible-data}
 \omega_0,\ \frac{\omega_0}{r}
 \in L^1(\R^3)\cap L^\infty(\R^3).
\end{equation}
Let \(u_0\) be the finite-energy Biot--Savart velocity generated by
\(\boldsymbol\Omega_0=\omega_0\mathbf{e_\theta}\); standard Biot--Savart and
Lorentz-space estimates show that \eqref{eq:admissible-data} implies
that \((u_0,\omega_0)\) satisfies the hypotheses of the global
well-posedness result in \cite{Danchin}.  We
denote the resulting unique global axisymmetric
solution without swirl by \((u,\omega)\); the regularity and transport
properties used below are recalled at the beginning of
Section~\ref{sec:prelim}.  A meridional scalar is identified with its
axisymmetric extension whenever an \(L^p(\R^3)\)-norm is used.

For two-point formulas, we abbreviate
\[
 X=(r,z),\qquad \bar X=(\bar r,\bar z),\qquad
 \dd X=\dd r\dd z,\qquad \dd\bar X=\dd\bar r\dd\bar z,
\]
and use, for example, \(\omega(X,t)=\omega(r,z,t)\).
Write
\begin{equation}\label{eq:S-def}
 S(X,\bar X):=
 \frac{(r-\bar r)^2+(z-\bar z)^2}{r\bar r},
 \qquad X,\bar X\in\Pi.
\end{equation}
Define
\begin{equation}\label{eq:F-def}
 F(s):=\int_0^\pi
 \frac{\cos\theta}{[2(1-\cos\theta)+s]^{1/2}}\,\dd\theta,
 \qquad s>0.
\end{equation}
We further define the auxiliary stream function
\begin{equation}\label{associated stream function}
    \psi(r,z,t):=\frac1{2\pi}\iint_\Pi
 \sqrt{r\bar r}\,F(S(X,\bar X))
 [-\omega(\bar r,\bar z,t)]\,\dd\bar r\dd\bar z.
\end{equation}
Then the scalar equation is
\begin{equation}\label{eq:conservation-law}
 \partial_t\omega+\partial_r(u^r\omega)+\partial_z(u^z\omega)=0,
 \qquad
 \partial_r u^r+\partial_z u^z=-\frac{u^r}{r},
\end{equation}
where \(u^r=r^{-1}\partial_z\psi\) and
\(u^z=-r^{-1}\partial_r\psi\).
Hence
\begin{align}
 u^r(r,z,t)
 &=\frac1{\pi r}\iint_\Pi
 \frac{z-\bar z}{\sqrt{r\bar r}}\,
 F'(S(X,\bar X))[-\omega(\bar r,\bar z,t)]
 \,\dd\bar r\dd\bar z,
 \label{eq:ur-full}\\
 u^z(r,z,t)
 &=-\frac1{4\pi r^2}\iint_\Pi\sqrt{r\bar r}\,
 F(S(X,\bar X))[-\omega(\bar r,\bar z,t)]
 \,\dd\bar r\dd\bar z\notag\\
 &\quad-\frac1{\pi r}\iint_\Pi\sqrt{r\bar r}
 \left(\frac{r-\bar r}{r\bar r}
       -\frac{S(X,\bar X)}{2r}\right)F'(S(X,\bar X))
 [-\omega(\bar r,\bar z,t)]\,\dd\bar r\dd\bar z.
 \label{eq:uz-full}
\end{align}
The symmetry and sign conditions in \eqref{eq:symmetry-sign} are
preserved by the flow.  Hence, for every \(t\geq0\),
\[
 \omega(r,-z,t)=-\omega(r,z,t),\qquad
 \omega(r,z,t)\leq0\quad(z>0).
\]
On \(\Pp\), set
\begin{equation}\label{eq:positive-density}
 \xi(r,z,t):=-\frac{\omega(r,z,t)}{r}\geq0,\qquad
 \xi_0:=\xi(\cdot,0).
\end{equation}
Thus \(-\omega=\partial_r u^z-\partial_z u^r\geq0\) on \(\Pp\), and
\(-\omega\) satisfies the same conservation law as \(\omega\).  We
retain the notation \(-\omega\) for this nonnegative density on
\(\Pp\).

Define the positive upper-half-plane moments
\begin{align}
 m_0&:=\iint_{\Pp}[-\omega(r,z,t)]\,\dd r\dd z,
 \label{eq:mass}\\
 P(t)&:=\iint_{\Pp}r^2[-\omega(r,z,t)]\,\dd r\dd z.
 \label{eq:P}
\end{align}
The mass \(m_0>0\) is conserved, and \(P(0)>0\).  Define the radial
support by
\begin{equation}\label{eq:radial-support-def}
 \mathcal R_\omega(t):=
 \inf\{R>0:\omega(r,z,t)=0\ \text{for }r>R,\ z>0\}.
\end{equation}
Here and below, support is understood modulo null sets.

If \(0\leq f\leq-\omega(\cdot,t)\) on \(\Pp\), let \(u_f\) be the
Biot--Savart velocity generated by the odd extension of
\(-f \mathbf{e_\theta}\), and write
\begin{equation}\label{eq:truncated-energy-intro}
 \Ecal(f):=\frac12\|u_f\|_{L^2(\R^3)}^2.
\end{equation}
In particular,
\(\Ecal(-\omega(\cdot,t))=\frac12\|u_0(\cdot)\|_{L^2(\R^3)}^2\).
The positive-kernel formula for \(\Ecal\) is given in
Section~\ref{sec:positive-kernels}.

\begin{theorem}\label{thm:main}
For compactly supported initial data $\omega_0(\cdot)\not\equiv0$ satisfying \eqref{eq:symmetry-sign} and \eqref{eq:admissible-data},
let \(\omega(\cdot,t)\) be the corresponding global solution of \eqref{eq:conservation-law}.  Then
\begin{equation}\label{eq:three-halves-main}
 \lim_{t\to\infty}
 \frac{P(t)[\log(2+t)]^{5/2}}{(1+t)^{3/2}}
 =+\infty.
\end{equation}
Consequently,
\begin{equation}\label{eq:three-halves-support}
 \lim_{t\to\infty}
 \frac{\mathcal R_\omega(t)[\log(2+t)]^{5/4}}
 {(1+t)^{3/4}}
 =+\infty.
\end{equation}
Moreover, for every \(0<\eta<1\), there exist
\(c_\eta>0\) and \(T_\eta>1\), depending only on \(\eta\) and the
initial data, such that, for every \(T\ge T_\eta\),
\[
 \left|
 \left\{t\in[T,2T]:
 \mathcal R_\omega(t)\ge c_\eta t
 \right\}
 \right|\ge(1-\eta)T.
\]
In addition, for every \(A>0\),
\begin{equation}\label{eq:typical-support}
 \lim_{T\to\infty}\frac1T
 \left|\left\{t\in[T,2T]:
 \mathcal R_\omega(t)\geq
 A\,t(\log t)^{-2}\right\}\right|=1.
\end{equation}
\end{theorem}

\begin{figure}[t]
\centering
\includegraphics[width=.94\textwidth]{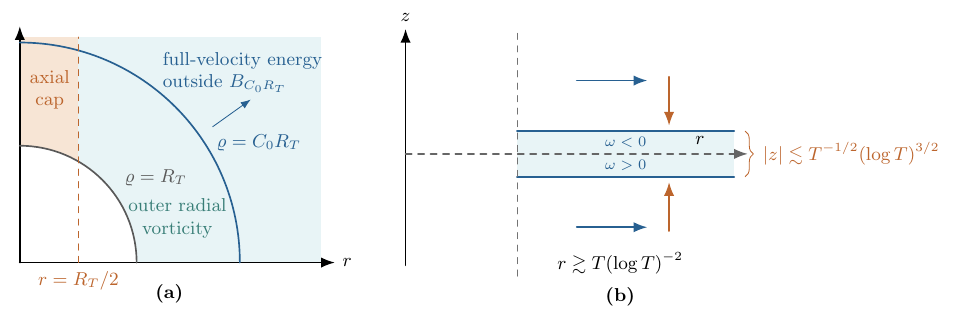}
\caption{The localization argument first transfers the exterior
velocity energy to vorticity at large cylindrical radius and then
places the energy-generating part in a shrinking neighborhood of the
symmetry plane.  In panel~(b), the two signs come from odd reflection
across \(z=0\).  The arrows describe radial escape and approach to
\(z=0\),
not tracked particle trajectories, and the shaded layer marks the
vorticity truncation rather than the spatial location of the velocity
energy density.}
\label{fig:energy-localization}
\end{figure}

The last conclusion of Theorem~\ref{thm:main} only detects the presence
of vorticity at large cylindrical radius.  A sharper use of the two
mixed moments and Lemma~\ref{lem:harmonic-tail} gives an energetic
version of the collision picture: for a density-one set of times, a
vorticity truncation generating an arbitrarily large fraction of the
conserved kinetic energy lies simultaneously at large \(r\) and in a
shrinking neighborhood of \(z=0\), as illustrated in
Figure~\ref{fig:energy-localization}(b).

\begin{theorem}\label{thm:escape-collision}
Let \(\omega(\cdot,t)\) be the global solution as in Theorem \ref{thm:main}.  For every
\(A>0\) and \(0<\delta<1\), there exists \(C_{A,\delta}>0\), depending
only on \(A\), \(\delta\), and the initial data, such that
\begin{equation}\label{eq:escape-collision}
\begin{aligned}
 \lim_{T\to\infty}\frac1T
 \Bigg|\Bigg\{t\in[T,2T]:\
 \Ecal\!\left(
 [-\omega(\cdot,t)]
 \ind_{\substack{
 D_{A,\delta}(T)
 }}
 \right)
 \geq \frac{1-\delta}{2}
 \|u_0(\cdot)\|_{L^2(\R^3)}^2
 \Bigg\}\Bigg|=1,
\end{aligned}
\end{equation}
where
\begin{equation*}
  D_{A,\delta}(T)=\{(r,z):  r>A T(\log T)^{-2},\,
 0<z<C_{A,\delta}T^{-1/2}(\log T)^{3/2}\}.
\end{equation*}
\end{theorem}

The velocity induced by the vorticity truncation in
\eqref{eq:escape-collision} carries an arbitrarily large fraction of the conserved kinetic energy.  Here \(\Ecal\) is a full-space kinetic energy, not the
integral of \(|u|^2/2\) over the displayed region.  The assertion is
instantaneous and Eulerian: it neither follows a fixed material set nor
proves persistence of a coherent dipole head, and it does not exclude
circulation being shed into a tail.

\begin{theorem}\label{thm:general-Lp}
Let \(\omega(\cdot,t)\) be the global solution fixed above.  Then
\begin{equation}\label{eq:general-Lp-full}
 \lim_{t\to\infty}
 \frac{[\log(2+t)]^{25/12}}{(1+t)^{1/4}}
 \inf_{1\leq p\leq\infty}
 \|\boldsymbol\Omega(t)\|_{L^p(\R^3)}
 =+\infty.
\end{equation}
Moreover, for every \(A>0\),
\begin{equation}\label{eq:general-Lp-typical}
 \lim_{T\to\infty}\frac1T
 \left|\left\{t\in[T,2T]:
 \inf_{1\leq p\leq\infty}
 \|\boldsymbol\Omega(t)\|_{L^p(\R^3)}
 \geq A t^{1/2}{(\log t)^{-3/2}} 
 \right\}\right|=1.
\end{equation}
\end{theorem}

We next record the additional bounds available for vortex patches.
Let \(D_0\subset\Pp\) be a bounded measurable set of positive measure
with \(\operatorname{dist}(D_0,\partial\Pp)>0\).  We also consider patch
data of the form
\begin{equation}\label{eq:patch-data}
 \xi_0(r,z)=\ind_{D_0}(r,z)\quad\text{in }\Pp,
\end{equation}
or, equivalently,
\(\omega_0(r,z)=-r\ind_{D_0}(r,z)\) on \(\Pp\), with odd extension across
\(z=0\).  The transported representative has the form
\begin{equation}\label{eq:patch-evolution}
 \xi(r,z,t)=\ind_{D_t}(r,z),\qquad
 D_t=\mathcal X_{t,0}(D_0),
\end{equation}
where \(\mathcal X_{t,0}\) is the meridional flow from time \(0\) to
time \(t\).  No regularity of
\(\partial D_0\) is assumed.

\begin{theorem}\label{thm:patch-Lp}
Let \(\omega(\cdot,t)\) be the global solution fixed above for the vortex patch initial
data \eqref{eq:patch-data}.  For every \(0<\eta<1\), there exist
\(c_\eta>0\) and \(T_\eta>1\), depending only on \(\eta\) and the
initial data, such that, for every \(T\ge T_\eta\),
\[
 \left|
 \left\{t\in[T,2T]:
 \mathcal R_\omega(t)
 =\|\boldsymbol\Omega(t)\|_{L^\infty(\R^3)}
 \ge c_\eta t
 \right\}
 \right|\ge(1-\eta)T.
\]
The vorticity maximum also satisfies the full-time estimate
\begin{equation}\label{eq:patch-Linfty-pointwise}
 \lim_{t\to\infty}
 \frac{\|\boldsymbol\Omega(t)\|_{L^\infty(\R^3)}
 [\log(2+t)]^{5/4}}{(1+t)^{3/4}}
 =+\infty.
\end{equation}
More generally,
\begin{subequations}\label{eq:patch-Lp-pointwise}
\begin{align}
 &\lim_{t\to\infty}\inf_{1\leq p<2}
 \frac{\|\boldsymbol\Omega(t)\|_{L^p(\R^3)}
 [\log(2+t)]^{5/(2p)}}
 {(1+t)^{(8p-7)/(6p)}}
 =+\infty,
 \label{eq:patch-Lp-pointwise-low}\\
 &\,\,\lim_{t\to\infty}\inf_{2\leq p\leq\infty}
 \frac{\|\boldsymbol\Omega(t)\|_{L^p(\R^3)}
 [\log(2+t)]^{5/4}}
 {(1+t)^{3/4}}
 =+\infty.\label{eq:patch-Lp-pointwise-high}
\end{align}
\end{subequations}

We also have the following family of density-one bounds.  We use the
convention \(1/\infty=0\).
For every fixed \(1\leq p\leq\infty\) and every \(A>0\),
\begin{equation}\label{eq:patch-Lp-typical}
 \lim_{T\to\infty}\frac1T
 \left|\left\{t\in[T,2T]:
 \|\boldsymbol\Omega(t)\|_{L^p(\R^3)}
 \geq A    t^{1-\frac1{2p}}{(\log t)^{\frac1{2p}-2} }   \right\}\right|=1.
\end{equation}
\end{theorem}

The general estimate \eqref{eq:general-Lp-full} also applies to
patches.  For \(1\leq p<2\), it and
\eqref{eq:patch-Lp-pointwise-low} give complementary bounds.

\subsection{Mixed moments and proof strategy}

Our proof is organized around a new pair of mixed moments.  Their role
is to turn conservation of kinetic energy into quantitative radial
escape, and then to convert that escape into growth of the radial
moment and the vorticity norms.  The argument consists of the following five steps.

\begin{itemize}

\item \textbf{Monotonicity of \(P(t)\).}
Since the relative vorticity \(\xi=-\omega/r\) is transported,
amplification of \(-\omega=r\xi\) can occur only through motion away
from the symmetry axis.  The monotonicity formula in
\cite[Theorem~1.1]{CJ} and the axis identity in
\cite[Section~4]{EY} give
\[
 P'(t)=\int_0^\infty r\,u^r(r,0,t)^2\,\dd r.
\]
Thus radial induction cannot cancel at the level of \(P\).  The main
difficulty is to show that the conserved kinetic energy forces a
substantial amount of vorticity to reach a large cylindrical radius.

\item\textbf{Two monotone mixed moments.}
The first weight is
\begin{equation}\label{eq:Phi}
 \Phi(r,z)=\int_0^z\sqrt{r^2+s^2}\,\dd s.
\end{equation}
Define
\begin{equation}\label{MPhi}
   M_\Phi(t):=\iint_{\Pi_+}\Phi(r,z)[-\omega(r,z,t)]\,\dd r\dd z.
\end{equation}
Proposition~\ref{prop:mixed} proves that
\[
 M_\Phi'(t)
 =-\frac12\int_0^\infty z|u^z(0,z,t)|^2\,\dd z
  -\iint_{\Pp}
  \frac{|zu^r-ru^z|^2}{r\sqrt{r^2+z^2}}\,\dd r\dd z
 \leq0.
\]
Since \(\Phi\geq rz\), it follows that
\begin{equation}\label{eq:strategy-rz}
 \iint_{\Pp}rz[-\omega(r,z,t)]\,\dd r\dd z\leq M_\Phi(0).
\end{equation}
Together with the energy estimates, this bound will exclude the
energy-generating part from lying at large \(z\) with \(r\) remaining
small.

The linear support estimate is not obtained by converting the lower
bound for \(P(t)\) into a support bound; that route loses a square root.
Instead, we apply \eqref{eq:multiplier} to a compactly supported cutoff
of \(R^{-1}z(r^2+z^2)\).  After integration on \([T,2T]\),
\eqref{eq:Phi-perfect-square} bounds the kinetic energy in
\(\{\varrho\leq R\}\) by
\[
 C_\varepsilon R+C_\varepsilon R
 [M_\Phi(T)-M_\Phi(2T)]+C\varepsilon T.
\]
First choosing \(\varepsilon\) small and then taking
\(R=a_\eta T\) with \(a_\eta>0\) small puts a fixed amount of velocity
energy outside \(B_{C_0R}\) on a proportion at least \(1-\eta\) of
\([T,2T]\).  Lemma~\ref{lem:harmonic-tail} then shows that vorticity
inside \(B_R\) cannot generate this exterior energy, and
\eqref{eq:strategy-rz} removes the high axial cap in
Figure~\ref{fig:energy-localization}(a).  Hence nonzero vorticity must
lie at \(r\gtrsim T\).  This direct use of the conserved kinetic energy
is the key step that gives linear, rather than square-root, radial
escape.

The second weight is chosen to control both velocity components with
the logarithmic loss needed below.  With
\(\varrho=(r^2+z^2)^{1/2}\), set
\[
 g(\varrho)=\frac{\varrho}{\operatorname{arsinh}\varrho},\qquad
 \Psi(r,z)=\lambda_0\Phi(r,z)+zg(\varrho),\qquad
 M_\Psi(t)=\iint_{\Pp}\Psi[-\omega]\,\dd r\dd z.
\]
Writing
\[
 u_{\rm rad}=\frac{ru^r+zu^z}{\varrho},\qquad
 u_{\rm tan}=\frac{-zu^r+ru^z}{\varrho},
\]
Proposition~\ref{prop:critical-mixed} shows that \(M_\Psi\) is
non-increasing and that \(\mathfrak D=-M_\Psi'\) satisfies
\[
 \mathfrak D(t)\gtrsim
 \iint_{\Pp}\left\{
 \frac{r}{\varrho[1+\operatorname{arsinh}\varrho]^2}|u_{\rm rad}|^2
 +\frac{\varrho}{r}|u_{\rm tan}|^2
 \right\}\,\dd r\dd z,\qquad
 \int_0^\infty\mathfrak D(t)\,\dd t\leq M_\Psi(0).
\]

\item\textbf{Kinetic energy outside expanding balls.}
Proposition~\ref{prop:critical-mixed} gives
\[
 \iint_{\{\varrho\leq R\}\cap\Pp}r|u|^2\,\dd r\dd z
 \lesssim R[\log(2+R)]^2\mathfrak D(t).
\]
Fix \(\kappa>0\) and take
\(R=R_T=\kappa T(\log T)^{-2}\), as in
\eqref{eq:critical-RT}.  The integrability of
\(\mathfrak D\) implies that, outside an exceptional set of
\(o_\kappa(T)\) times in \([T,2T]\), a fixed portion of the conserved
kinetic energy lies outside \(B_{C_0R}\), where \(C_0>1\) is a fixed
absolute constant.

It remains to identify the vorticity responsible for that energy.
The velocity generated by vorticity supported in \(B_R\) is a harmonic
gradient in \(\R^3\setminus B_R\).  Its monopole vanishes, and a
spherical-harmonic expansion yields
\[
 \|u_{\rm in}\|_{L^2(\R^3\setminus B_{C_0R})}
 \leq C_0^{-3/2}\|u_{\rm in}\|_{L^2(\R^3)}.
\]
Choosing \(C_0\) large shows that the velocity generated by the outer
vorticity has a fixed amount of energy.  Keeping track of the small
constants makes this amount arbitrarily close to the total kinetic
energy.  Estimate
\eqref{eq:strategy-rz} first removes the part lying at large \(z\) but
small \(r\).  A second use of the same mixed moment, together with the
logarithmic kernel estimate, removes the portion above
\[
 z\sim T^{-1/2}(\log T)^{3/2}.
\]
Consequently, for every \(\delta>0\), at least \(1-\delta\) of the
kinetic energy is generated by the vorticity truncation in
\eqref{eq:escape-collision} for all but
\(o_{\kappa,\delta}(T)\) times in \([T,2T]\).

\item\textbf{Radial escape and concentration near the symmetry plane.}
The preceding estimates give
Theorem~\ref{thm:escape-collision}.  Using the compactly supported
multiplier above together with Lemma~\ref{lem:harmonic-tail} instead
gives linear radial escape on any prescribed fixed proportion of
\([T,2T]\), which proves the linear support assertion in
Theorem~\ref{thm:main}.  For patches, the support radius is exactly the
vorticity maximum, giving the linear conclusion in
Theorem~\ref{thm:patch-Lp}.

\item\textbf{Full-time growth.}
The estimates in \cite[Sections~3.1--3.3]{GMT}, combined with
the quantitative logarithmic estimate in
Lemma~\ref{lem:log-Lp}, give
\[
 P'(t)\gtrsim \kappa^{1/2}T^{1/2}(\log T)^{-5/2}
\]
at every time outside the exceptional set just described.  Integrating
over \([T,2T]\) and then using the
arbitrariness of \(\kappa\) proves \eqref{eq:three-halves-main}.  This
converts the preceding escape estimate into a full-time pointwise
lower bound for \(P\).

The full-time support bound \eqref{eq:three-halves-support} follows immediately from the moment estimate.  For
general data, combining it with the radial-velocity estimate in
\cite[Proposition~A.1]{LimJeong} and preservation of the
vorticity-support volume gives \eqref{eq:general-Lp-full}.  For
patches, the identity \(-\omega=r\) on the transported patch gives the
estimates in Theorem~\ref{thm:patch-Lp}.
\end{itemize}

The transport statements, including those for patches with rough
boundary, are justified using the velocity estimates in
\cite{Danchin} and the renormalized transport theory of
DiPerna--Lions \cite{DiPernaLions}.  The axis regularizations and
cutoff limits needed for the two mixed moments are given in
Appendix~\ref{app:mixed-justification}.

\noindent\textbf{Organization of the paper:}
Section~\ref{sec:prelim} fixes the measure conventions and collects the
transport properties, moment identities, and
kernel and energy estimates used in the proof.
Section~\ref{sec:proof} establishes the two mixed-moment identities,
proves the linear support estimate, and shows that the vorticity
truncation in \eqref{eq:escape-collision} generates almost all of the
conserved kinetic energy.  It then derives the radial-moment and support
bounds and concludes with the general and patch
\(L^p\)-estimates.
Appendix~\ref{app:mixed-justification} supplies the axis
regularizations and cutoff limits for the mixed moments.

\section{Preliminaries}
\label{sec:prelim}

We will apply several results and calculations from Choi--Jeong \cite{CJ},
Egamberganov--Yao \cite{EY}, and Gustafson--Miller--Tsai \cite{GMT}.
For completeness, and to fix the conventions and normalizations needed
later, we include the relevant details.

\subsection{Transport and moments}

Unless stated otherwise, \(L^p(\R^3)\)-norms of axisymmetric scalars use
three-dimensional measure:
\begin{equation}\label{eq:3d-axisymmetric-measure}
 \|h\|_{L^p(\R^3)}^p
 =2\pi\iint_\Pi |h(r,z)|^p r\,\dd r\dd z,
 \qquad 1\leq p<\infty.
\end{equation}
For \(p=\infty\), the norm is the essential supremum on \(\R^3\).
Integrals over \(\Pi\) or \(\Pp\), including \(m_0\) and \(P\), are
unweighted.

On every finite time interval, the velocity is log-Lipschitz, locally
uniformly in time, by \cite{Danchin}.  The axisymmetric velocity
estimates in \cite{Danchin} also give
\begin{equation}\label{eq:danchin-axis-control}
 \sup_{0\leq t\leq T}
 \left\|\frac{u^r(t)}r\right\|_{L^\infty(\R^3)}
 <\infty\qquad(T<\infty).
\end{equation}
Since \(\omega(t)\in L^1(\R^3)\cap L^\infty(\R^3)\), the Biot--Savart law gives
\(u(t)\in W^{1,p}(\R^3)\) for every finite \(p\geq2\).  We use the
continuous Biot--Savart representatives on \(r=0\) and \(z=0\); they
agree with the local Sobolev traces.

For completeness, we record why the kinetic energy is conserved in
this class, including for vortex patches with rough boundary.  On every finite
interval \([0,T]\), the preceding estimates imply
\[
 u,\nabla u\in L^\infty(0,T;L^3(\R^3)).
\]
Let \(\rho_\varepsilon\) be a spatial mollifier and set
\[
 \mathcal C_\varepsilon
 :=(u\otimes u)_\varepsilon
   -u_\varepsilon\otimes u_\varepsilon.
\]
Mollifying the Euler equation and taking the
\(L^2(\R^3)\)-inner product with \(u_\varepsilon\) gives
\begin{equation}\label{eq:mollified-energy-balance}
 \frac12\frac{\dd}{\dd t}
 \|u_\varepsilon(t)\|_{L^2(\R^3)}^2
 =\int_{\R^3}\nabla u_\varepsilon:
 \mathcal C_\varepsilon\,\dd x.
\end{equation}
The pressure and transport terms vanish because
\(\nabla\cdot u_\varepsilon=0\).  Moreover,
\[
 \|\nabla u_\varepsilon(t)\|_{L^3(\R^3)}
 \leq\|\nabla u(t)\|_{L^3(\R^3)},\qquad
 \|\mathcal C_\varepsilon(t)\|_{L^{3/2}(\R^3)}
 \longrightarrow0
\]
for almost every \(t\); the second convergence follows from the strong
continuity of translations in \(L^3(\R^3)\).  The product of these two
norms is uniformly integrable on \([0,T]\).  Integrating
\eqref{eq:mollified-energy-balance}, passing to the limit by dominated
convergence, and using
\(u_\varepsilon(t)\to u(t)\) in \(L^2(\R^3)\) proves
\[
 \|u(t)\|_{L^2(\R^3)}=\|u(0)\|_{L^2(\R^3)}
 \qquad(t\geq0).
\]

The renormalized transport theory \cite{DiPernaLions} and the
log-Lipschitz flow imply
\begin{equation}\label{eq:relative-vorticity-bound}
 0\leq\xi(r,z,t)\leq\|\xi_0\|_{L^\infty(\Pp)},
 \qquad
 0\leq-\omega(r,z,t)\leq
 r\|\xi_0\|_{L^\infty(\Pp)}
 \quad\text{in }\Pp.
\end{equation}
The flow preserves \(r\,\dd r\dd z\), leaves \(\Pp\) invariant, and
keeps the vorticity compactly supported on bounded time intervals.
Oddness in \(z\) gives
\[
 u^r(r,-z,t)=u^r(r,z,t),\qquad
 u^z(r,-z,t)=-u^z(r,z,t),
\]
and hence \(u^z(r,0,t)=0\).

If \(\phi\) is locally Lipschitz near the vorticity support, then
\(t\mapsto\iint_{\Pp}\phi[-\omega]\,\dd X\) is locally absolutely continuous
and
\begin{equation}\label{eq:moment-derivative}
 \frac{\dd}{\dd t}\iint_{\Pp}\phi(X)[-\omega(X,t)]\,\dd X
 =\iint_{\Pp}\nabla\phi(X)\cdot u(X,t)[-\omega(X,t)]\,\dd X
 \quad\text{for almost every }t.
\end{equation}
All moment functions below are taken in their locally absolutely continuous
representatives; their differential identities therefore hold for almost
every \(t\).
Compact support and the Biot--Savart law also give, uniformly on bounded
time intervals,
\begin{equation}\label{eq:far-field-decay}
 |u(x,t)|=O(|x|^{-2}),\qquad
 |\nabla u(x,t)|=O(|x|^{-3}),\qquad
 |u^r(r,z,t)|=O(r|x|^{-3})
 \quad (|x|\to\infty).
\end{equation}

Because \(\omega(\cdot,t)\) is odd in \(z\), the full-plane function
\(-\omega(\cdot,t)\) is the odd extension of its nonnegative restriction
to \(\Pp\).  For $F$ defined in \eqref{eq:F-def}, we use
\cite[Lemma~3.1]{GMT}:
\begin{equation}\label{eq:Fprime-two-sided}
 \frac{c}{s(1+s)^{3/2}}
 \leq-F'(s)\leq\frac{C}{s(1+s)^{3/2}},
 \qquad s>0,
\end{equation}
and, with the same constant \(c\),
\[
 -F'(s)\geq\frac{c}{6s},\qquad 0<s\leq2,
\]
because \((1+s)^{3/2}\leq3^{3/2}<6\) on this interval.

The kinetic energy is conserved.  Using symmetry in \(z\) and integration in
the angular variable, we write it as
\begin{equation}\label{eq:E0}
 E_0:=\iint_{\Pp}r|u(r,z,t)|^2\,\dd r\dd z
 =\frac1{4\pi}\|u(t)\|_{L^2(\R^3)}^2.
\end{equation}
Indeed,
\[
 \|u(t)\|_{L^2(\R^3)}^2
 =2\pi\iint_\Pi r|u(r,z,t)|^2\,\dd r\dd z
 =4\pi\iint_{\Pp}r|u(r,z,t)|^2\,\dd r\dd z,
\]
because \(|u|^2\) is even in \(z\).
Since \(\omega_0\not\equiv0\), one has \(E_0>0\).
From this point onward, \(c\) and \(C\) denote positive constants that
may depend on the fixed initial data, in particular on
\(m_0,\|\xi_0\|_{L^\infty(\Pp)},P(0),E_0,M_\Phi(0)\), and \(M_\Psi(0)\), but not on
\(t,T,R\), or \(\kappa\).  Dependence on an additional parameter,
such as \(p\), is indicated by a subscript.

Choi--Jeong \cite{CJ} proved the monotonicity of \(P\) through a
positive double-integral representation of \(P'\).
Egamberganov--Yao \cite[Section~4.1]{EY} gave a different proof based
on the axis identity
\begin{equation}
 P'(t)=2\iint_{\Pp}
 r u^r(r,z,t)[-\omega(r,z,t)]\,\dd r\dd z
 =\int_0^\infty r\,u^r(r,0,t)^2\,\dd r.
 \label{eq:Pprime-axis}
\end{equation}
The first equality in \eqref{eq:Pprime-axis}, together with
\eqref{eq:relative-vorticity-bound} and Cauchy--Schwarz, gives
\[
\begin{aligned}
 P'(t)
 &\leq2\iint_{\Pp}
 r|u^r(r,z,t)|[-\omega(r,z,t)]\,\dd r\dd z\\
 &\leq2
 \left(\iint_{\Pp}r|u^r(r,z,t)|^2\,\dd r\dd z\right)^{1/2}
 \left(\iint_{\Pp}r[-\omega(r,z,t)]^2\,\dd r\dd z\right)^{1/2}\\
 &\leq2\bigl(\|\xi_0\|_{L^\infty(\Pp)}E_0P(t)\bigr)^{1/2}.
\end{aligned}
\]
In deducing the last inequality, we have used
\[
 \iint_{\Pp}r|u^r|^2\,\dd r\dd z\leq E_0,
 \qquad
 \iint_{\Pp}r[-\omega]^2\,\dd r\dd z
 \leq\|\xi_0\|_{L^\infty(\Pp)}P(t),
\]
where the second estimate follows pointwise from
\(-\omega\leq r\|\xi_0\|_{L^\infty(\Pp)}\).  

\noindent Since \(P(t)\geq P(0)>0\),
division by \(2P(t)^{1/2}\) gives
\[
 \frac{\dd}{\dd t}P(t)^{1/2}
 \leq\bigl(\|\xi_0\|_{L^\infty(\Pp)}E_0\bigr)^{1/2}
\]
for almost every \(t\).  Integrating from \(0\) to \(t\) yields
\begin{equation}\label{eq:P-upper}
 P(t)\leq
 \bigl(P(0)^{1/2}
 +(\|\xi_0\|_{L^\infty(\Pp)}E_0)^{1/2}t\bigr)^2
 \leq C(1+t)^2.
\end{equation}
The preceding argument and \eqref{eq:P-upper} are proved in
\cite[Section~3]{EY}.

Splitting \eqref{eq:ur-full} into the two half-planes gives
\begin{equation}\label{eq:ur-axis-explicit}
 u^r(r,0,t)=\frac2\pi\iint_{\Pp}
 \frac{\bar z}{r\sqrt{r\bar r}}
 \left[-F'\!\left(
 \frac{(r-\bar r)^2+\bar z^2}{r\bar r}
 \right)\right][-\omega(\bar r,\bar z,t)]\,\dd\bar r\dd\bar z.
\end{equation}
Indeed, at \(z=0\) the kernel argument is unchanged by
\(\bar z\mapsto-\bar z\).  Under this reflection, both
\(-\bar z\) in \eqref{eq:ur-full} and the odd function
\(-\omega(\bar r,\bar z,t)\) change sign.  The contribution of the lower
half-plane therefore equals that of the upper half-plane, which accounts
for the factor two in \eqref{eq:ur-axis-explicit}.
\subsection{Positive kernels and vorticity truncations}
\label{sec:positive-kernels}

\begin{figure}[H]
\centering
\includegraphics[width=.82\textwidth]{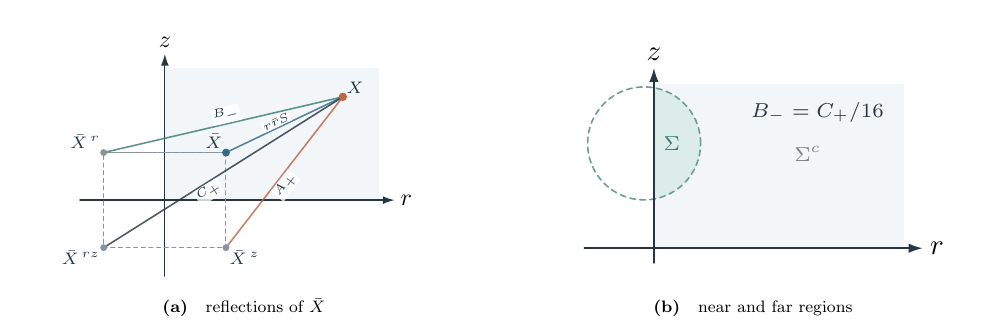}
\caption{Reflected-point geometry for the positive kernels.
Panel~(b) illustrates the near--far decomposition introduced below.}
\label{fig:kernel-reflections}
\end{figure}

For \(X=(r,z)\), \(\bar X=(\bar r,\bar z)\in\Pp\), define
\begin{align}
 L(X,\bar X)&:=\log(2+S(X,\bar X)^{-1}),\label{eq:S-L}\\
 A_+(X,\bar X)&:=(r-\bar r)^2+(z+\bar z)^2,
 &B_-(X,\bar X)&:=(r+\bar r)^2+(z-\bar z)^2,\notag\\
 C_+(X,\bar X)&:=(r+\bar r)^2+(z+\bar z)^2,\notag\\
 K(X,\bar X)&:=
 \frac{(z+\bar z)(r\bar r)^2}
 {A_+(X,\bar X)C_+(X,\bar X)^{3/2}},
 &e(X,\bar X)&:=
 \frac{z\bar z(r\bar r)^2L(X,\bar X)}
 {A_+(X,\bar X)B_-(X,\bar X)^{3/2}}.
 \label{eq:K-e}
\end{align}
Set also
\[
 \bar S(X,\bar X):=\frac{A_+(X,\bar X)}{r\bar r}.
\]
Let
\(\bar X^{\,z}=(\bar r,-\bar z)\),
\(\bar X^{\,r}=(-\bar r,\bar z)\), and
\(\bar X^{\,rz}=(-\bar r,-\bar z)\).  Then
\[
 r\bar rS=|X-\bar X|^2,\qquad
 A_+=|X-\bar X^{\,z}|^2,\qquad
 B_-=|X-\bar X^{\,r}|^2,\qquad
 C_+=|X-\bar X^{\,rz}|^2;
\]
see Figure~\ref{fig:kernel-reflections}(a).

Let \(f\) be any bounded, compactly supported, nonnegative measurable
function on \(\Pp\).  Denote by \(f^{\rm odd}\) its odd extension in \(z\).
As in \eqref{eq:truncated-energy-intro}, let \(u_f\) be the
Biot--Savart velocity generated by \(-f^{\rm odd}\mathbf{e_\theta}\), so that
\(\Ecal(f)=\frac12\|u_f\|_{L^2(\R^3)}^2\).
The vector field \(-f^{\rm odd}\mathbf{e_\theta}\) is divergence-free in the sense
of distributions and belongs to \(L^1(\R^3)\cap L^2(\R^3)\), so
\(u_f\in L^2(\R^3)\).
For smooth \(f\), the constant and the reflected difference in the
energy formula can be checked directly.  The generated scalar vorticity
is \(-f^{\rm odd}\); let \(\psi_f\) be the associated stream function in
\eqref{associated stream function}.  Angular integration and integration by parts give
\[
\begin{aligned}
 \Ecal(f)
 &=\frac12\|u_f\|_{L^2(\R^3)}^2
 =\pi\iint_\Pi r|u_f|^2\,\dd X
 =\pi\iint_\Pi\frac{|\nabla\psi_f|^2}{r}\,\dd X\\
 &=-\pi\iint_\Pi\psi_f
 \left\{\partial_r\!\left(\frac{\partial_r\psi_f}{r}\right)
       +\partial_z\!\left(\frac{\partial_z\psi_f}{r}\right)\right\}
 \,\dd X\\
 &=\pi\iint_\Pi\psi_f f^{\rm odd}\,\dd X
 =\frac12\iint_{\Pi^2}\sqrt{r\bar r}\,
 F(S(X,\bar X))f^{\rm odd}(X)f^{\rm odd}(\bar X)
 \,\dd X\dd\bar X.
\end{aligned}
\]
The four choices of the signs of \(z\) and \(\bar z\) contribute two
copies of the same-side kernel \(F(S)\) and two copies, with the
opposite sign, of the reflected kernel \(F(\bar S)\).  The factor
\(1/2\) in the full-plane formula therefore gives
\begin{equation}\label{eq:energy-static-exact}
 \Ecal(f)=
 \iint_{\Pp^2}\sqrt{r\bar r}\,
 [F(S(X,\bar X))-F(\bar S(X,\bar X))]
 f(X)f(\bar X)\,\dd X\dd\bar X.
\end{equation}
For measurable \(f\), choose nonnegative functions
\(f_n\in C_c^\infty(\Pp)\) by first truncating away from \(r=0\),
\(z=0\), and spatial infinity and then mollifying inside \(\Pp\).
Their odd axisymmetric extensions satisfy
\[
 f_n^{\rm odd}\longrightarrow f^{\rm odd}
 \quad\text{in }L^1(\R^3)\cap L^2(\R^3),
\]
and hence in \(L^{6/5}(\R^3)\).
The Hardy--Littlewood--Sobolev inequality gives convergence of the
corresponding Biot--Savart velocities in \(L^2(\R^3)\) and convergence
of the Newtonian energy pairings.  Thus
\eqref{eq:energy-static-exact} applies to every measurable truncation
\(0\leq f\leq-\omega(\cdot,t)\), without any regularity assumption on
the truncating set.

The point of the two kernels in \eqref{eq:K-e} is now visible:
\(K\) measures a positive contribution to the stretching rate
\(P'(t)\), whereas \(e\) controls the kinetic energy generated by a
vorticity truncation.

We next derive both estimates, including the reflected denominators.
They correspond to \cite[(3.8) and (3.9)]{GMT}; writing them out here
also fixes every normalization that enters the near--far argument.
Recall that \(S\) and \(\bar S\) are the original and \(z\)-reflected
kernel arguments.  Directly from their definitions,
\begin{equation}\label{eq:denominator-normalizations}
 A_+=r\bar r\,\bar S,\qquad
 B_-=r\bar r(S+4),\qquad
 C_+=r\bar r(\bar S+4),\qquad
 \bar S-S=\frac{4z\bar z}{r\bar r}.
\end{equation}

First split the full-plane formula \eqref{eq:ur-full} into
\(\bar z>0\) and \(\bar z<0\).  In the lower half-plane, reflect
\(\bar z\) and use the oddness of \(-\omega\).  This gives
\begin{equation}\label{eq:ur-half-plane}
 u^r(r,z,t)=\frac1{\pi r}\iint_{\Pp}
 \frac{(z-\bar z)F'(S(X,\bar X))
       -(z+\bar z)F'(\bar S(X,\bar X))}
 {\sqrt{r\bar r}}[-\omega(\bar X,t)]\,\dd\bar X.
\end{equation}
Substitution into the first equality in \eqref{eq:Pprime-axis} yields
\begin{align}
 P'(t)
 ={}&\frac2\pi\iint_{\Pp^2}
 \frac{(z-\bar z)F'(S(X,\bar X))}{\sqrt{r\bar r}}\,
 [-\omega(X,t)][-\omega(\bar X,t)]\,\dd X\dd\bar X\notag\\
 &-\frac2\pi\iint_{\Pp^2}
 \frac{(z+\bar z)F'(\bar S(X,\bar X))}{\sqrt{r\bar r}}\,
 [-\omega(X,t)][-\omega(\bar X,t)]\,\dd X\dd\bar X.
 \label{eq:Pprime-before-antisymmetry}
\end{align}
The first integrand is antisymmetric under
\(X\leftrightarrow\bar X\), and thus its integral is zero.  For completeness,
this exchange is legitimate: if \(d=|X-\bar X|\), then
\eqref{eq:Fprime-two-sided} gives
\[
 \frac{|z-\bar z|}{\sqrt{r\bar r}}
 |F'(S(X,\bar X))|
 \leq C\frac{|z-\bar z|(r\bar r)^2}
 {d^2(r\bar r+d^2)^{3/2}}
 \leq C\frac{\sqrt{r\bar r}}{d}.
\]
At a fixed time the vorticity is bounded and compactly supported, and
\(d^{-1}\) is locally integrable in the two-dimensional meridional
variables.  Therefore
\begin{equation}\label{eq:Pprime-reflected-Fprime}
 P'(t)=\frac2\pi\iint_{\Pp^2}
 \frac{z+\bar z}{\sqrt{r\bar r}}\,
 [-F'(\bar S(X,\bar X))]
 [-\omega(X,t)][-\omega(\bar X,t)]\,\dd X\dd\bar X.
\end{equation}
The lower bound estimate in \eqref{eq:Fprime-two-sided} and
\eqref{eq:denominator-normalizations} give, one line at a time,
\[
\begin{aligned}
 \frac1{\sqrt{r\bar r}}[-F'(\bar S)]
 &\geq
 \frac{c}{\sqrt{r\bar r}\,\bar S(1+\bar S)^{3/2}}\\
 &=c\frac{(r\bar r)^2}
 {A_+(A_++r\bar r)^{3/2}}
 \geq c\frac{(r\bar r)^2}{A_+C_+^{3/2}},
\end{aligned}
\]
because \(A_++r\bar r\leq A_++4r\bar r=C_+\).  Inserting this bound in
\eqref{eq:Pprime-reflected-Fprime} proves
\begin{equation}
 P'(t)\geq c\iint_{\Pp^2}
 K(X,\bar X)[-\omega(X,t)][-\omega(\bar X,t)]
 \,\dd X\dd\bar X.
 \label{eq:Pprime-K}
\end{equation}

We now turn to the energy kernel.  The only one-dimensional estimate
needed is
\begin{equation}\label{eq:one-dimensional-energy-kernel}
 \int_s^{s+\delta}\frac{\dd\tau}
 {\tau(1+\tau)^{3/2}}
 \leq C\frac{\delta}{s+\delta}
 \frac{\log(2+s^{-1})}{(s+4)^{3/2}},
 \qquad s,\delta>0.
\end{equation}
To show \eqref{eq:one-dimensional-energy-kernel}, we distinguish the three elementary cases.  If \(0<\delta\leq s\), the
left-hand side is at most
\(\delta/[s(1+s)^{3/2}]\); use
\(s+\delta\leq2s\), \(s+4\leq4(1+s)\), and
\(\log(2+s^{-1})\geq\log2\).  If \(\delta>s\) and \(0<s\leq1\), then
\(\delta/(s+\delta)\geq1/2\), \(s+4\sim1\), and
\[
 \int_s^{s+\delta}\frac{\dd\tau}
 {\tau(1+\tau)^{3/2}}
 \leq\int_s^1\frac{\dd\tau}{\tau}
     +\int_1^\infty\tau^{-5/2}\,\dd\tau
 \leq C\log(2+s^{-1}).
\]
Finally, if \(\delta>s\) and \(s\geq1\), then \(s+4\sim s\) and
\[
 \int_s^{s+\delta}\frac{\dd\tau}
 {\tau(1+\tau)^{3/2}}
 \leq\int_s^\infty\tau^{-5/2}\,\dd\tau
 =\frac23s^{-3/2}.
\]
This proves \eqref{eq:one-dimensional-energy-kernel}.

Take \(s=S(X,\bar X)\) and
\(\delta=\bar S(X,\bar X)-S(X,\bar X)
=4z\bar z/(r\bar r)\).  Since \(F'<0\), the upper bound estimate in
\eqref{eq:Fprime-two-sided} and
\eqref{eq:one-dimensional-energy-kernel} give
\[
\begin{aligned}
 F(S)-F(\bar S)
 =\int_S^{\bar S}[-F'(\tau)]\,\dd\tau
 \leq C\frac{\bar S-S}{\bar S}
 \frac{\log(2+S^{-1})}{(S+4)^{3/2}}.
\end{aligned}
\]
Multiplying by \(\sqrt{r\bar r}\) and using
\eqref{eq:denominator-normalizations} in each factor, we obtain
\begin{align}
 \sqrt{r\bar r}\,[F(S)-F(\bar S)]
 &\leq C\sqrt{r\bar r}\,
 \frac{4z\bar z/(r\bar r)}{A_+/(r\bar r)}
 \frac{L(X,\bar X)}
 {[B_-/(r\bar r)]^{3/2}}\notag\\
 &=C\frac{z\bar z\,(r\bar r)^2L(X,\bar X)}
 {A_+B_-^{3/2}}
 =Ce(X,\bar X).
 \label{eq:energy-kernel-pointwise}
\end{align}
Substituting this pointwise estimate in
\eqref{eq:energy-static-exact} gives, for every measurable
\(0\leq\eta\leq-\omega(\cdot,t)\),
\begin{equation}
 \Ecal(\eta)\leq C\iint_{\Pp^2}
 e(X,\bar X)\eta(X)\eta(\bar X)\,\dd X\dd\bar X.
 \label{eq:energy-e}
\end{equation}
The constant $C$ in \eqref{eq:energy-kernel-pointwise} is independent of
\(\eta\); hence so is the constant in \eqref{eq:energy-e}.

We will also use the following crude estimate:
\begin{equation}\label{eq:e-crude}
 e(X,\bar X)\leq C\sqrt{r\bar r}\,L(X,\bar X).
\end{equation}
Indeed, \(A_+\geq4z\bar z\) and \(B_-\geq4r\bar r\) give
\[
 \frac{z\bar z(r\bar r)^2}{A_+B_-^{3/2}}
 \leq\frac14\,\frac{(r\bar r)^2}{(4r\bar r)^{3/2}}
 =\frac1{32}\sqrt{r\bar r}.
\]

We next write out the near--far comparison used below; it is the
three-dimensional form of the decomposition in
\cite[Section~3.3]{GMT}.  We fix the separation ratio \(1/4\) once
and for all and define
\begin{equation}\label{eq:Sigma}
 \Sigma:=\{(X,\bar X)\in\Pp^2:
 B_-(X,\bar X)\leq C_+(X,\bar X)/16\}.
\end{equation}
For fixed \(\bar X\), this is the intersection of \(\Pp\) with an
Apollonius disk.  We call \(\Sigma\) the near region and
\(\Sigma^c\) the far region; see
Figure~\ref{fig:kernel-reflections}(b).

On \(\Sigma\),
\[
 C_+-B_-=4z\bar z,\qquad
 \frac{15}{16}C_+\leq4z\bar z,
\]
and hence
\[
 C_+^{1/2}\leq
 \frac8{\sqrt{15}}(z\bar z)^{1/2}
 \leq C(z\bar z)^{1/2}.
\]
Using also \(A_+\geq4z\bar z\),
\((r\bar r)^2\leq B_-^2/16\), and
\(B_-^{1/2}\leq C_+^{1/2}/4\), we obtain
\[
\begin{aligned}
 e(X,\bar X)
 &=\frac{z\bar z(r\bar r)^2L}{A_+B_-^{3/2}}
 \leq\frac{(r\bar r)^2L}{4B_-^{3/2}}
 \leq CB_-^{1/2}L\\
 &\leq C\,C_+^{1/2}L
 \leq C(z\bar z)^{1/2}L.
\end{aligned}
\]
Thus
\begin{equation}\label{eq:near-GMT}
 e(X,\bar X)\leq C(z\bar z)^{1/2}L(X,\bar X).
\end{equation}
Next, on \(\Sigma^c\), for every \(0<\nu<1\), we have
\begin{equation}\label{eq:far-GMT}
 e(X,\bar X)\leq C(r\bar r)^{(1-\nu)/2}
 (z\bar z)^{\nu/2}K(X,\bar X)^\nu L(X,\bar X).
\end{equation}
The constant in \eqref{eq:far-GMT} may be chosen uniformly
for \(0<\nu<1\).  Indeed, on \(\Sigma^c\) the inequality
\(B_-^{3/2}>C_+^{3/2}/64\) reduces the estimate to the factor
\[
 \frac{z\bar z(r\bar r)^2}{A_+C_+^{3/2}}.
\]
The inequalities \(A_+\geq4z\bar z\), \(C_+\geq4r\bar r\), and
\(z+\bar z\geq2(z\bar z)^{1/2}\) give
\[
 \frac{z\bar z(r\bar r)^2}{A_+C_+^{3/2}}
 \leq\frac1{32}(r\bar r)^{1/2},
 \qquad
 \frac{z\bar z(r\bar r)^2}{A_+C_+^{3/2}}
 =\frac{z\bar z}{z+\bar z}K
 \leq\frac12(z\bar z)^{1/2}K.
\]
Interpolating these two bounds proves \eqref{eq:far-GMT}, with a constant
uniform for \(0<\nu<1\).  In the application below, \(\nu\) will depend
on the dyadic time scale and lie in \((1/2,2/3)\).

\subsection{Logarithmic kernel estimate}

The logarithmic-kernel estimate in \cite[Lemma~3.3]{GMT} treats fixed
\(p\).  That form is
not sufficient here, because the H\"older exponent varies with the
dyadic time scale.  We therefore sharpen their estimate slightly by
making its dependence on \(p\) explicit.  Recall from
\eqref{eq:S-L} that
\(L(X,\bar X)=\log(2+S(X,\bar X)^{-1})\), and write
\[
 \|L\|_{L^p(\eta\otimes\eta)}
 :=
 \left(\iint_{\Pp^2}|L(X,\bar X)|^p
 \eta(X)\eta(\bar X)\,\dd X\dd\bar X\right)^{1/p}.
\]

\begin{lemma}\label{lem:log-Lp}
There is
\(C_L=C(m_0,\|\xi_0\|_{L^\infty(\Pp)})\) such that, for every real \(p\geq1\),
every \(t\geq0\), and every
\(0\leq\eta(r,z)\leq-\omega(r,z,t)\) on \(\Pp\),
\begin{equation}\label{eq:log-Lp-quantitative}
 \|L\|_{L^p(\eta\otimes\eta)}
 \leq C_L\bigl[p+\log(2+P(t))\bigr].
\end{equation}
Consequently, for fixed \(p<\infty\),
\begin{equation}\label{eq:log-Lp}
 \|L\|_{L^p(\eta\otimes\eta)}
 \leq C_{p,L}\log(2+P(t)).
\end{equation}
\end{lemma}

\begin{proof}
Let
\[
 \dd\mu_\eta=\eta(X)\eta(\bar X)\,\dd X\dd\bar X,\qquad
 d=|X-\bar X|.
\]
Then \(\mu_\eta(\Pp^2)\leq m_0^2\).  Throughout this proof, \(C\)
depends only on \(m_0\) and
\(\|\xi_0\|_{L^\infty(\Pp)}\).

For every \(\ell\geq\log4\),
\begin{equation}\label{eq:L-level-implies-S}
 L>\ell\quad\Longrightarrow\quad
 d^2\leq2e^{-\ell}r\bar r.
\end{equation}
Indeed,
\[
 L>\ell
 \quad\Longrightarrow\quad
 S^{-1}>e^\ell-2\geq\frac12e^\ell,
\]
which yields \eqref{eq:L-level-implies-S}, since
\(d^2=Sr\bar r\).

Fix \(R_*\geq1\).  Chebyshev's inequality for the radial moment gives
\begin{equation}\label{eq:L-tail-large-radius}
 \mu_\eta\{L>\ell,\ \max(r,\bar r)>R_*\}
 \leq2m_0P(t)R_*^{-2}.
\end{equation}
Indeed,
\[
 \iint_{\{r>R_*\}}\eta(X)\,\dd X
 \leq R_*^{-2}\iint_{\Pp}r^2\eta(X)\,\dd X
 \leq R_*^{-2}P(t).
\]
Applying this estimate to \(r\) and \(\bar r\) gives
\eqref{eq:L-tail-large-radius}.

On \(r,\bar r\leq R_*\), for each fixed \(X\),
\eqref{eq:L-level-implies-S} confines \(\bar X\) to a planar ball of
radius \(CR_*e^{-\ell/2}\).  Since
\(\eta(\bar X)\leq
\|\xi_0\|_{L^\infty(\Pp)}\bar r
\leq\|\xi_0\|_{L^\infty(\Pp)}R_*\),
\begin{equation}\label{eq:L-tail-bounded-radius}
 \mu_\eta\{L>\ell,\ r,\bar r\leq R_*\}
 \leq C\|\xi_0\|_{L^\infty(\Pp)}m_0R_*^3e^{-\ell}.
\end{equation}

Choose
\[
 R_*=(2+P(t))^{1/5}e^{\ell/5}.
\]
This is at least one.  Combining
\eqref{eq:L-tail-large-radius} and
\eqref{eq:L-tail-bounded-radius} yields
\begin{equation}\label{eq:L-exponential-tail}
 \mu_\eta\{L>\ell\}
 \leq C(2+P(t))^{3/5}e^{-2\ell/5}.
\end{equation}

Now set
\[
 \ell_0=C_{\rm tail}+\frac32\log(2+P(t)),
\]
where
\(C_{\rm tail}=C_{\rm tail}
(m_0,\|\xi_0\|_{L^\infty(\Pp)})\)
is sufficiently large.  Then \(\ell_0\geq\log4\), and
\eqref{eq:L-exponential-tail} gives
\[
 \mu_\eta\{L>\ell_0+s\}\leq m_0^2e^{-2s/5},\qquad s\geq0.
\]
For \(Y=(L-\ell_0)_+\), the layer-cake formula gives
\[
 \|Y\|_{L^p(\mu_\eta)}^p
 \leq m_0^2p\int_0^\infty s^{p-1}e^{-2s/5}\,\dd s
 =m_0^2\Gamma(p+1)(5/2)^p.
\]
Since \(L\leq\ell_0+Y\) and
\(\mu_\eta(\Pp^2)\leq m_0^2\), Stirling's estimate and
Minkowski's inequality imply
\[
 \|L\|_{L^p(\eta\otimes\eta)}
 \leq C_L\bigl[p+\log(2+P(t))\bigr],
\]
where
\(C_L=C_L(m_0,\|\xi_0\|_{L^\infty(\Pp)})\).
This proves \eqref{eq:log-Lp-quantitative}, and \eqref{eq:log-Lp} follows
for fixed \(p\) from
\(\log(2+P(t))\geq\log2\).
\end{proof}

\section{Proof of the main results}\label{sec:proof}

We first establish the two identities for the mixed moments and combine
them with \eqref{eq:multiplier}, applied to a compactly supported cutoff,
to transfer the conserved kinetic energy to vorticity at large radius.
The estimates in Section~\ref{sec:positive-kernels} then yield growth
of \(P(t)\).  Finally, we convert the moment and localization bounds into
the \(L^p\)-estimates for general data and vortex patches.

\subsection{The multiplier identity and the first mixed moment}
\label{sec:mixed}

We look for a multiplier \(\phi\) for which the quadratic form in
\eqref{eq:multiplier} is a non-positive square.  This requirement leads to
\(\phi_z=(r^2+z^2)^{1/2}\).  Together with \(\phi(r,0)=0\), this gives
the function \(\Phi\) in \eqref{eq:Phi}.  We first establish the multiplier
identity and let \(\phi=\Phi\) to obtain a uniform bound for
\(\iint_{\Pp}rz[-\omega(r,z,t)]\,\dd r\dd z\).  We then use the same
identity with the logarithmically weighted multiplier
\(\lambda_0\Phi+z\varrho/\operatorname{arsinh}\varrho\), which controls both
spherical components of the velocity.

\begin{lemma}\label{lem:multiplier}
For every \(\phi\in C_c^2([0,\infty)^2)\), define
\[
 \mathcal M_\phi(t):=
 \iint_{\Pp}\phi(r,z)[-\omega(r,z,t)]\,\dd r\dd z.
\]
Then \(\mathcal M_\phi\) is locally absolutely continuous and
\begin{align}
 \mathcal M_\phi'(t)
 ={}&\frac12\int_0^\infty
 \phi_r(r,0)u^r(r,0,t)^2\,\dd r
 -\frac12\int_0^\infty
 \phi_z(0,z)u^z(0,z,t)^2\,\dd z
 \notag\\
 &+\iint_{\Pp}
 \left(\phi_{rz}(r,z)-\frac{\phi_z(r,z)}{r}\right)
 u^r(r,z,t)^2\,\dd r\dd z
 -\iint_{\Pp}
 \phi_{rz}(r,z)u^z(r,z,t)^2\,\dd r\dd z
 \notag\\
 &+\iint_{\Pp}
 [\phi_{zz}(r,z)-\phi_{rr}(r,z)]
 u^r(r,z,t)u^z(r,z,t)\,\dd r\dd z\notag\\
 &+\iint_{\Pp}\frac{\phi_r(r,z)}{r}
 u^r(r,z,t)u^z(r,z,t)\,\dd r\dd z.
 \label{eq:multiplier}
\end{align}
\end{lemma}

\begin{proof}
By \eqref{eq:moment-derivative},
\begin{equation}\label{eq:multiplier-first-step}
 \mathcal M_\phi'(t)=\iint_{\Pp}
 (\phi_ru^r+\phi_zu^z)
 (\partial_ru^z-\partial_zu^r)\,\dd r\dd z.
\end{equation}
The two integrals in \eqref{eq:multiplier-first-step} involving \(u^r\)
and its derivatives, in which the derivative falls on the same velocity
component, give
\[
\begin{aligned}
 -\iint_{\Pp}\phi_ru^r\partial_zu^r\,\dd r\dd z
 &=\frac12\int_0^\infty
   \phi_r(r,0)u^r(r,0,t)^2\,\dd r
   +\frac12\iint_{\Pp}\phi_{rz}(u^r)^2\,\dd r\dd z,\\
 \iint_{\Pp}\phi_zu^z\partial_ru^z\,\dd r\dd z
 &=-\frac12\int_0^\infty
   \phi_z(0,z)u^z(0,z,t)^2\,\dd z
   -\frac12\iint_{\Pp}\phi_{rz}(u^z)^2\,\dd r\dd z.
\end{aligned}
\]
The boundary terms at infinity vanish because \(\phi\) is compactly
supported.

For the integrals in \eqref{eq:multiplier-first-step} involving products
of \(u^r\) and \(u^z\) and their derivatives, integration by parts in
\(r\) and \(z\), together with \(u^r(0,z,t)=0\) and
\(u^z(r,0,t)=0\), gives
\begin{align}
 &\iint_{\Pp}
 \bigl(\phi_ru^r\partial_ru^z-\phi_zu^z\partial_zu^r\bigr)
 \,\dd r\dd z\notag\\
 &\quad=\iint_{\Pp}(\phi_{zz}-\phi_{rr})u^ru^z\,\dd r\dd z
 -\iint_{\Pp}\phi_r(\partial_ru^r)u^z\,\dd r\dd z
 +\iint_{\Pp}\phi_z(\partial_zu^z)u^r\,\dd r\dd z.
 \label{eq:multiplier-cross-first}
\end{align}
Using
\[
 \partial_ru^r+\partial_zu^z=-\frac{u^r}{r},
\]
the last two terms respectively become
\[
 \begin{aligned}
 -\iint_{\Pp}\phi_r(\partial_ru^r)u^z\,\dd r\dd z
 &=-\frac12\iint_{\Pp}\phi_{rz}(u^z)^2\,\dd r\dd z
   +\iint_{\Pp}\frac{\phi_r}{r}u^ru^z\,\dd r\dd z,\\
 \iint_{\Pp}\phi_z(\partial_zu^z)u^r\,\dd r\dd z
 &=\frac12\iint_{\Pp}\phi_{rz}(u^r)^2\,\dd r\dd z
   -\iint_{\Pp}\frac{\phi_z}{r}(u^r)^2\,\dd r\dd z.
\end{aligned}
\]
Adding these identities gives \eqref{eq:multiplier}.  The regularity
stated in Section~\ref{sec:prelim}, together with compactly supported
cutoffs, justifies the integrations by parts and the passage to the two
axes.  Local absolute continuity follows from
\eqref{eq:moment-derivative}.
\end{proof}

\begin{proposition}\label{prop:mixed}
The mixed moment \(M_\Phi(t)\) defined in \eqref{MPhi} is non-increasing.  In particular,
\begin{equation}\label{eq:rz-bound}
 \iint_{\Pp}rz[-\omega(r,z,t)]\,\dd r\dd z\leq M_\Phi(0)
 \qquad(t\geq0).
\end{equation}
\end{proposition}

\begin{proof}
Put \(\varrho=(r^2+z^2)^{1/2}\).  Direct differentiation away from the
axis gives
\[
 \Phi_z=\varrho,\qquad
 \Phi_r=r\,\operatorname{arsinh}(z/r),\qquad
 \Phi_{rz}=\frac r\varrho,\qquad
 \Phi_{zz}-\Phi_{rr}+\frac{\Phi_r}{r}=\frac{2z}{\varrho}.
\]
The regularization at the axis and the removal of the outer cutoff are
justified in Lemma~\ref{lem:Phi-justification} of
Appendix~\ref{app:mixed-justification}.  Thus, although \(\Phi\) is not
compactly supported, \eqref{eq:multiplier} remains valid with
\(\phi=\Phi\).  Substituting the identities above into
\eqref{eq:multiplier} gives
\begin{align}
 &\frac{\dd}{\dd t}\iint_{\Pp}
 \Phi(r,z)[-\omega(r,z,t)]\,\dd r\dd z\notag\\
 &\quad=-\frac12\int_0^\infty z\,u^z(0,z,t)^2\,\dd z
 -\iint_{\Pp}
 \frac{[zu^r(r,z,t)-ru^z(r,z,t)]^2}
 {r\sqrt{r^2+z^2}}\,\dd r\dd z
 \leq0.
 \label{eq:Phi-perfect-square}
\end{align}
Indeed, \(\Phi_r(r,0)=0\), \(\Phi_z(0,z)=z\), and the bulk quadratic form is
\[
 -\frac{z^2}{r\varrho}u^r(r,z,t)^2
 -\frac r\varrho u^z(r,z,t)^2
 {}+\frac{2z}{\varrho}u^r(r,z,t)u^z(r,z,t)
 =-\frac{[zu^r(r,z,t)-ru^z(r,z,t)]^2}{r\varrho}.
\]
  Hence
\eqref{eq:Phi-perfect-square} holds for the solution considered here.
It follows that \(M_\Phi\) is non-increasing, and
\(\Phi(r,z)\geq rz\) gives \eqref{eq:rz-bound}.
\end{proof}

\subsection{A second mixed moment controlling both velocity components}
\label{sec:critical-moment}

The square in \eqref{eq:Phi-perfect-square} controls
\(zu^r-ru^z\), namely the component of the velocity tangent to the
centered circles \(\varrho=\mathrm{const}\) in the meridional plane.
It gives no control of the radial component.  To localize the full
kinetic energy, we therefore introduce a second weight whose
dissipation controls both components, at the cost of a logarithmic
loss.  For the remainder of the proof,
\(\varrho=(r^2+z^2)^{1/2}\).  Let
\begin{equation}\label{eq:critical-g}
 g(\varrho):=
 \begin{cases}
  \varrho/\operatorname{arsinh}\varrho,&\varrho>0,\\
  1,&\varrho=0.
 \end{cases}
\end{equation}
Since
\(g(\varrho)=1+\varrho^2/6+O(\varrho^4)\) as
\(\varrho\to0\), the function \(zg(\varrho)\) is smooth at the
origin.

\begin{proposition}\label{prop:critical-mixed}
There is an absolute constant \(\lambda_0>0\) such that, for
\begin{equation}\label{eq:critical-Psi}
 \Psi(r,z):=\lambda_0\Phi(r,z)+zg(\varrho),\qquad
 M_\Psi(t):=\iint_{\Pp}\Psi(r,z)[-\omega(r,z,t)]\,\dd r\dd z,
\end{equation}
the function \(M_\Psi\) is locally absolutely continuous, nonnegative,
and non-increasing.  More precisely, set
\begin{equation}\label{eq:spherical-components}
 u_{\rm rad}:=\frac{ru^r+zu^z}{\varrho},\qquad
 u_{\rm tan}:=\frac{-zu^r+ru^z}{\varrho}.
\end{equation}
Then there exists a universal constant \(c_0>0\) such that, for almost
every \(t>0\),
\begin{equation}\label{eq:critical-coercivity}
 \mathfrak D(t):=-M_\Psi'(t)
 \geq c_0\iint_{\Pp}
 \left\{
  \frac{r}{\varrho[1+\operatorname{arsinh}\varrho]^2}
  |u_{\rm rad}|^2
  +\frac{\varrho}{r}|u_{\rm tan}|^2
 \right\}\,\dd r\dd z.
\end{equation}
Consequently,
\begin{equation}\label{eq:critical-budget}
 \int_0^\infty\mathfrak D(t)\,\dd t\leq M_\Psi(0).
\end{equation}
\end{proposition}

\begin{proof}
For a multiplier \(\phi\), denote by \(Q_\phi(u^r,u^z)\) the sum of the
four volume integrals on the right-hand side of
\eqref{eq:multiplier}.  Proposition~\ref{prop:mixed} gives
\begin{equation}\label{eq:Q-Phi-spherical}
 Q_\Phi=-\frac{\varrho}{r}u_{\rm tan}^2.
\end{equation}
For the form generated by \(zg(\varrho)\), write
\begin{equation}\label{eq:w-dg-def}
 \ell:=\operatorname{arsinh}\varrho,\qquad
 w:=\frac g\varrho-g',\qquad
 d_g:=\frac{2g}\varrho-\varrho g''.
\end{equation}
Since \(g(\varrho)=\varrho/\ell\), direct differentiation gives
\begin{equation}
 w(\varrho)=\frac{\varrho}
 {\sqrt{1+\varrho^2}\,\ell^2}
 =\frac{\tanh\ell}{\ell^2}>0.
 \label{eq:g-derivatives}
\end{equation}
Set
\[
 \gamma=\frac r\varrho,\qquad \sigma=\frac z\varrho.
\]
Similarly, substitution of \(\phi=zg(\varrho)\) into
\eqref{eq:multiplier}, followed by the orthogonal change of variables
\eqref{eq:spherical-components}, gives
\begin{equation}\label{eq:Qg-spherical}
 Q_{zg}
 =-\gamma w u_{\rm rad}^2
  -\frac{g'+\sigma^2g/\varrho}{\gamma}u_{\rm tan}^2
  +\sigma d_g u_{\rm rad}u_{\rm tan}.
\end{equation}
Thus the negative of the bulk form generated by \(\Psi\) is represented
in the variables \((u_{\rm rad},u_{\rm tan})\) by
\begin{equation}\label{eq:critical-matrix}
 \mathcal H=
 \begin{pmatrix}
  \gamma w&-\sigma d_g/2\\[1mm]
  -\sigma d_g/2&
  \dfrac{\lambda_0+g'+\sigma^2g/\varrho}{\gamma}
 \end{pmatrix}.
\end{equation}

It remains to choose \(\lambda_0\).  The quantity that measures the
possible loss in the determinant is
\begin{equation}\label{eq:Theta-def}
 \Theta(\varrho):=
 \frac{d_g(\varrho)^2}{4w(\varrho)}
 -g'(\varrho)-\frac{g(\varrho)}{\varrho}.
\end{equation}
The Taylor expansion of \(\operatorname{arsinh}\varrho\) at zero gives
\begin{equation}\label{eq:Theta-small}
 w=\frac1\varrho-\frac{\varrho}{6}+O(\varrho^3),
 \qquad d_g=\frac2\varrho+O(\varrho^3),\qquad
 \Theta=-\frac{\varrho}{3}+O(\varrho^3).
\end{equation}
Thus the apparent singularity of \(\Theta\) at the origin is
removable.

At infinity, write \(\varrho=\sinh\ell\), so that
\(\sqrt{1+\varrho^2}=\cosh\ell\).
Using \eqref{eq:w-dg-def}--\eqref{eq:g-derivatives} yields the exact
forms
\[
 w=\frac{\tanh\ell}{\ell^2},\qquad
 d_g=\frac2\ell
 +\frac{\tanh\ell(1+\operatorname{sech}^2\ell)}{\ell^2}
 -\frac{2\tanh^2\ell}{\ell^3}.
\]
Since \(\tanh\ell=1+O(e^{-2\ell})\), these formulas yield
\begin{equation}\label{eq:Theta-large}
 w=\frac{\tanh\ell}{\ell^2},\qquad
 d_g=\frac2\ell+\frac1{\ell^2}-\frac2{\ell^3}
   +O(e^{-2\ell}\ell^{-2}),
 \qquad
 \Theta=1-\frac1\ell+O(\ell^{-2}).
\end{equation}
Together with smoothness away from the origin,
\eqref{eq:Theta-small}--\eqref{eq:Theta-large} show that
\(\Theta\) is bounded above on \((0,\infty)\).  Moreover,
\begin{equation}\label{eq:gprime-positive}
 g'(\varrho)=\frac1\ell-\frac{\tanh\ell}{\ell^2},
 \qquad 0<g'(\varrho)<\frac1\ell,
\end{equation}
because \(0<\tanh\ell<\ell\).

Choose
\begin{equation}\label{eq:lambda0-choice}
 \lambda_0\geq2+\sup_{\varrho>0}\Theta(\varrho).
\end{equation}
By \eqref{eq:Theta-large},
\(\sup_{\varrho>0}\Theta(\varrho)\geq1\); hence
\eqref{eq:lambda0-choice} and \eqref{eq:gprime-positive} give
\[
 \lambda_0+g'(\varrho)\geq3,\qquad
 \lambda_0-\Theta(\varrho)\geq2.
\]
The determinant of \(\mathcal H\) can then be written as
\begin{align}
 \det\mathcal H
 &=w\left(\lambda_0+g'+\sigma^2\frac g\varrho\right)
   -\frac{\sigma^2d_g^2}{4}\notag\\
 &=w\left[(1-\sigma^2)(\lambda_0+g')
          +\sigma^2(\lambda_0-\Theta)\right]\geq w.
 \label{eq:critical-det}
\end{align}
Since \(w>0\), the upper-left entry and the determinant are positive,
so \(\mathcal H\) is positive definite.  Since
\(g/\varrho=1/\ell\) and \(0<g'<1/\ell\),
\begin{equation}\label{eq:critical-d-upper}
 \mathcal H_{22}\leq\frac{C(1+\ell^{-1})}{\gamma}.
\end{equation}
Completing the square in the two possible orders and averaging the
resulting lower bounds gives
\begin{align}
 (u_{\rm rad},u_{\rm tan})\mathcal H
 (u_{\rm rad},u_{\rm tan})^{\mathsf T}
 \geq\frac12\left(
 \frac{\det\mathcal H}{\mathcal H_{22}}u_{\rm rad}^2
 +\frac{\det\mathcal H}{\mathcal H_{11}}u_{\rm tan}^2\right).
 \label{eq:two-square-bound}
\end{align}
Since \(\mathcal H_{11}=\gamma w\), using
\eqref{eq:critical-det}--\eqref{eq:critical-d-upper} in
\eqref{eq:two-square-bound} gives
\begin{equation}\label{eq:matrix-coercivity-pre}
 (u_{\rm rad},u_{\rm tan})\mathcal H
 (u_{\rm rad},u_{\rm tan})^{\mathsf T}
 \geq c\left\{
  \gamma\frac{w}{1+\ell^{-1}}u_{\rm rad}^2
  +\frac1\gamma u_{\rm tan}^2
 \right\}.
\end{equation}
Finally,
\begin{equation}\label{eq:w-critical-weight}
 \frac{w}{1+\ell^{-1}}
 =\frac{\tanh\ell}{\ell(\ell+1)}
 \geq\frac{c}{(1+\ell)^2}.
\end{equation}
This proves the bulk estimate in \eqref{eq:critical-coercivity}.

The boundary terms also have the required sign:
\(\Psi_r(r,0)=0\) and
\begin{equation}\label{eq:critical-axis-coefficient}
 \Psi_z(0,z)=\lambda_0z+g(z)+zg'(z)>0.
\end{equation}
Lemma~\ref{lem:Psi-justification} in
Appendix~\ref{app:mixed-justification} verifies local integrability at
the axis and the cutoff limit at infinity.  It therefore gives
\begin{align}
 M_\Psi(t_2)-M_\Psi(t_1)
 ={}&-\int_{t_1}^{t_2}\iint_{\Pp}
 (u_{\rm rad},u_{\rm tan})\mathcal H
 (u_{\rm rad},u_{\rm tan})^{\mathsf T}
 \,\dd r\dd z\dd t\notag\\
 &-\frac12\int_{t_1}^{t_2}\int_0^\infty
 [\lambda_0z+g(z)+zg'(z)]|u^z(0,z,t)|^2\,\dd z\dd t.
 \label{eq:critical-integrated-identity}
\end{align}
Then \eqref{eq:matrix-coercivity-pre}--\eqref{eq:w-critical-weight}
prove \eqref{eq:critical-coercivity}.  Since \(\Psi\geq0\), integration
in time proves \eqref{eq:critical-budget}.
\end{proof}

Estimate \eqref{eq:critical-coercivity} controls the full kinetic energy
in every centered ball.  Together with the finite global $L^1$ bound on
$\mathfrak D$ in \eqref{eq:critical-budget}, it forces most of the
conserved energy to lie outside balls of almost linear radius for most
times in large dyadic intervals.

\begin{proposition}\label{prop:critical-inner-energy}
For every \(R\geq2\) and almost every \(t>0\),
\begin{equation}\label{eq:critical-inner-energy}
 \iint_{\Pp\cap\{\varrho\leq R\}}
 r|u(r,z,t)|^2\,\dd r\dd z
 \leq CR[\log(2+R)]^2\mathfrak D(t).
\end{equation}
\end{proposition}

\begin{proof}
Since \(|u|^2=u_{\rm rad}^2+u_{\rm tan}^2\) and
\(\varrho\leq R\), we have
\begin{align*}
 r u_{\rm rad}^2
 &\leq R[1+\operatorname{arsinh}R]^2
 \frac{r}
 {\varrho[1+\operatorname{arsinh}\varrho]^2}
 u_{\rm rad}^2,\\
 r u_{\rm tan}^2&\leq R\frac{\varrho}{r}u_{\rm tan}^2.
\end{align*}
Now use Proposition~\ref{prop:critical-mixed} and
 \(1+\operatorname{arsinh}R\leq C\log(2+R)\).
\end{proof}

\subsection{From total kinetic energy to exterior-vorticity energy}
\label{sec:critical-localization}

Proposition~\ref{prop:critical-inner-energy} controls the
\(L^2\)-energy of the full velocity \(u(t)\) inside a ball.  The
positive-kernel argument in the next subsection, however, requires an
\(L^2\)-lower bound for the velocity generated by the vorticity outside
a smaller ball.  The following exterior estimate allows us to pass
from the first statement to the second.

\begin{lemma}\label{lem:harmonic-tail}
Let \(\boldsymbol\Omega\in L^1(\R^3)\cap L^2(\R^3)\) be compactly
supported in \(B_R\), assume that
\(\nabla\cdot\boldsymbol\Omega=0\) in the sense of distributions, and set
\[
 v=\nabla\times(-\Delta)^{-1}\boldsymbol\Omega.
\]
Then, for every \(C>1\),
\begin{equation}\label{eq:harmonic-tail}
 \|v\|_{L^2(\R^3\setminus B_{CR})}
 \leq C^{-3/2}\|v\|_{L^2(\R^3)}.
\end{equation}
\end{lemma}

\begin{proof}
Since
\(\nabla\cdot v=0\) and
\(\nabla\times v=\boldsymbol\Omega\), one has
\[
 \nabla\cdot v=\nabla\times v=0
 \qquad\text{in }\R^3\setminus B_R.
\]
The exterior of \(B_R\) is simply connected.  Hence there exists a
harmonic function \(h\) on \(\R^3\setminus B_R\), normalized by
\(h(x)\to0\) as \(|x|\to\infty\), such that \(v=\nabla h\).

For every \(s>R\),
\[
 \int_{\partial B_s}v\cdot n\,\dd S
 =\int_{B_s}\nabla\cdot v\,\dd x=0.
\]
Let \(\{Y_{nm}\}\) be an \(L^2(\mathbb S^2)\)-orthonormal basis of
spherical harmonics.  The exterior expansion of \(h\) initially has
the form
\[
 h(x)=\sum_{n\geq0}\sum_m
 a_{nm}|x|^{-n-1}Y_{nm}(x/|x|).
\]
The flux identity forces the coefficient of the \(n=0\) mode to
vanish.  Therefore
\begin{equation}\label{eq:harmonic-expansion}
 h(x)=\sum_{n\geq1}\sum_m
 a_{nm}|x|^{-n-1}Y_{nm}(x/|x|).
\end{equation}
For every \(s>R\), this series and its gradient converge in the
exterior Dirichlet space on \(\{|x|>s\}\).
For \(s>R\), orthogonality gives
\[
 \int_{|x|>s}
 \left|\nabla\!\left(|x|^{-n-1}Y_{nm}(x/|x|)\right)\right|^2
 \,\dd x
 =(n+1)s^{-2n-1},
\]
and hence the explicit identity
\[
 \int_{|x|>s}|\nabla h|^2\,\dd x
 =\sum_{n\geq1}\sum_m
 (n+1)|a_{nm}|^2s^{-2n-1}.
\]
Consequently,
\[
\begin{aligned}
 \|v\|_{L^2(\R^3\setminus B_{Cs})}^2
 &=\sum_{n\geq1}\sum_m
 C^{-2n-1}(n+1)|a_{nm}|^2s^{-2n-1}\\
 &\leq C^{-3}
 \sum_{n\geq1}\sum_m
 (n+1)|a_{nm}|^2s^{-2n-1}\\
 &\leq C^{-3}\|v\|_{L^2(\R^3)}^2.
\end{aligned}
\]
Letting \(s\to R\) and taking square roots proves
\eqref{eq:harmonic-tail}.
\end{proof}

At times when \(\mathfrak D(t)\) is small,
Proposition~\ref{prop:critical-inner-energy} forces most of the
conserved kinetic energy outside a large ball.  The next result shows
that, for a density-one set of times, a fixed part of this energy is
generated by vorticity lying both outside a smaller ball and away from
the symmetry axis.  The second restriction is essential because
purely axial escape produces no vortex stretching.

\begin{proposition}\label{prop:critical-good-times}
There exists \(E_*>0\), depending only on the initial data, with the
following property.  For every fixed
\(\kappa>0\), there is \(T_{\rm loc}(\kappa)>1\) such that, for every
\(T\geq T_{\rm loc}(\kappa)\), setting
\begin{equation}\label{eq:critical-RT}
 R_T:=\kappa T(\log T)^{-2},
\end{equation}
there is a measurable set \(G_T\subset[T,2T]\) satisfying
\begin{equation}\label{eq:critical-good-measure}
 |[T,2T]\setminus G_T|=o_\kappa(T)
 \quad\text{as }T\to\infty,\qquad
 |G_T|\geq\frac{3T}{4},
\end{equation}
such that, for every \(t\in G_T\),
\begin{equation}\label{eq:critical-outer-energy}
 \Ecal\!\left(
 [-\omega(\cdot,t)]
 \ind_{\{\varrho>R_T,\ r>R_T/2\}}
 \right)\geq E_*.
\end{equation}
\end{proposition}

\begin{proof}
Recall that
\(\Ecal(f)=\frac12\|u_f\|_{L^2(\R^3)}^2\), where \(u_f\) is generated
by \(-f^{\rm odd}\mathbf{e_\theta}\).  Choose \(C_0>1\) once and for all
so large that \(C_0^{-3/2}\leq1/8\).
The proof has three steps, illustrated
in Figure~\ref{fig:energy-localization}(a): we first push the energy of
the full velocity outside \(B_{C_0R_T}\), then transfer it to
vorticity outside \(B_{R_T}\), and finally remove the high axial cap.

\smallskip
\noindent\emph{Step 1: Push full-velocity energy outside
\(B_{C_0R_T}\).}
Fix \(\kappa>0\).
By Proposition~\ref{prop:critical-inner-energy}, after increasing
\(T_{\rm loc}(\kappa)\) so that
\(2\leq C_0R_T\leq C_0T\) and
\(\log(2+C_0R_T)\leq \log T\),
\begin{equation}\label{eq:small-D-inner-energy}
 \mathfrak D(t)\leq
 \frac{E_0}{4C C_0R_T(\log T)^2}
 \quad\Longrightarrow\quad
 \iint_{\Pp\cap\{\varrho\leq C_0R_T\}}
 r|u(r,z,t)|^2\,\dd r\dd z
 \leq\frac{E_0}{4}.
\end{equation}
We call a time bad if the dissipation exceeds the threshold in
\eqref{eq:small-D-inner-energy}; thus the bad-time set is
\[
 B_T:=\left\{t\in[T,2T]:
 \mathfrak D(t)>
 \frac{E_0}{4C C_0R_T(\log T)^2}\right\}.
\]
Therefore
\[
 |B_T|\frac{E_0}{4C C_0R_T(\log T)^2}
 \leq\int_{B_T}\mathfrak D(t)\,\dd t
 \leq\int_T^{2T}\mathfrak D(t)\,\dd t,
\]
and hence
\begin{equation}\label{eq:critical-bad-measure}
 |B_T|\leq\frac{4C C_0R_T(\log T)^2}{E_0}
 \int_T^{2T}\mathfrak D(t)\,\dd t
 =C\kappa T
 \int_T^{2T}\mathfrak D(t)\,\dd t
 =o_\kappa(T).
\end{equation}
For fixed \(\kappa\), the last equality follows because
\(\mathfrak D\in L^1(0,\infty)\) implies
\(\int_T^{2T}\mathfrak D(t)\,\dd t\to0\) as \(T\to\infty\).

Now, we define the good-time set by
\[
 G_T:=[T,2T]\setminus B_T.
\]
  In particular,
\[
 |[T,2T]\setminus G_T|=|B_T|=o_\kappa(T),
\]
and every \(t\in G_T\) satisfies the hypothesis on the left of
\eqref{eq:small-D-inner-energy}.
Increase \(T_{\rm loc}(\kappa)\), if necessary, so that
\[
 |[T,2T]\setminus G_T|\leq T/4
\]
for every \(T\geq T_{\rm loc}(\kappa)\).  The symmetry in \(z\), angular
integration, \eqref{eq:E0}, and
\eqref{eq:small-D-inner-energy} imply, for \(t\in G_T\),
\begin{equation}\label{eq:full-velocity-outer}
 \|u(t)\|_{L^2(\R^3\setminus B_{C_0R_T})}
 \geq\frac{\sqrt3}{2}\|u(t)\|_{L^2(\R^3)}.
\end{equation}
Indeed,
\[
 \|u(t)\|_{L^2(B_{C_0R_T})}^2
 =4\pi\iint_{\Pp\cap\{\varrho\leq C_0R_T\}}
 r|u(r,z,t)|^2\,\dd r\dd z
 \leq\pi E_0
 =\frac14\|u(t)\|_{L^2(\R^3)}^2.
\]

\smallskip
\noindent\emph{Step 2: Transfer to vorticity outside \(B_{R_T}\).}
At a good time, \eqref{eq:full-velocity-outer} places most of the
full-velocity energy outside \(B_{C_0R_T}\).  We now show that this
energy cannot be generated mainly by the vorticity inside \(B_{R_T}\).
We split
\begin{equation*}
 -\omega(\cdot,t)=\eta_{\rm in}+\eta_{\rm out},\qquad
 \eta_{\rm in}:=[-\omega(\cdot,t)]\ind_{\{\varrho\leq R_T\}},\qquad
 \eta_{\rm out}:=[-\omega(\cdot,t)]\ind_{\{\varrho>R_T\}}.
\end{equation*}
Let \(u_{\rm in}\) and \(u_{\rm out}\) be their Biot--Savart velocities.
Writing \(\eta_{\rm in}^{\rm odd}\) for the odd extension in \(z\), the
vector field \(-\eta_{\rm in}^{\rm odd}\mathbf{e_\theta}\) is divergence-free
in the sense of distributions.  Hence Lemma~\ref{lem:harmonic-tail} gives
\begin{equation}\label{eq:inner-harmonic-tail}
 \|u_{\rm in}\|_{L^2(\R^3\setminus B_{C_0R_T})}
 \leq C_0^{-3/2}\|u_{\rm in}\|_{L^2(\R^3)}.
\end{equation}
Since \(u(t)=u_{\rm in}+u_{\rm out}\), the triangle inequality on the
exterior domain and \eqref{eq:full-velocity-outer} give
\[
\begin{aligned}
 \frac{\sqrt3}{2}\|u(t)\|_{L^2(\R^3)}
 &\leq
 \|u(t)\|_{L^2(\R^3\setminus B_{C_0R_T})}\\
 &\leq
 \|u_{\rm in}\|_{L^2(\R^3\setminus B_{C_0R_T})}
 +\|u_{\rm out}\|_{L^2(\R^3)}\\
 &\leq
 C_0^{-3/2}\|u_{\rm in}\|_{L^2(\R^3)}
 +\|u_{\rm out}\|_{L^2(\R^3)}\\
 &\leq
 C_0^{-3/2}\|u(t)\|_{L^2(\R^3)}
 +(1+C_0^{-3/2})\|u_{\rm out}\|_{L^2(\R^3)}.
\end{aligned}
\]
Our choice of \(C_0\) therefore implies
\(\|u_{\rm out}\|_{L^2(\R^3)}
\geq c\|u(t)\|_{L^2(\R^3)}\), and hence
\begin{equation}\label{eq:outer-vorticity-energy}
 \Ecal(\eta_{\rm out})
 =\frac12\|u_{\rm out}\|_{L^2(\R^3)}^2
 \geq cE_0.
\end{equation}
Here and for the remainder of this proof, \(c>0\) denotes a fixed
absolute constant for which \eqref{eq:outer-vorticity-energy} holds.

\smallskip
\noindent\emph{Step 3: Remove the high axial cap.}
We have shown that the velocity generated by the vorticity outside
\(B_{R_T}\) has a fixed amount of energy.  It remains to discard the
portion close to the axis.  Define
\begin{equation}\label{eq:high-cap}
 \eta_H:=[-\omega(\cdot,t)]
 \ind_{\{\varrho>R_T,\ r\leq R_T/2\}},\qquad
 \eta_{\rm rad}:=[-\omega(\cdot,t)]
 \ind_{\{\varrho>R_T,\ r>R_T/2\}}.
\end{equation}
Then \(\eta_{\rm out}=\eta_H+\eta_{\rm rad}\).
On \(\operatorname{supp}\eta_H\),
\(z\geq\sqrt3R_T/2\).  The estimates
\eqref{eq:energy-e} and \eqref{eq:e-crude}, followed by
Cauchy--Schwarz, give
\begin{equation}\label{eq:cap-energy-first}
 \Ecal(\eta_H)
 \leq C\left(\iint_{\Pp}r\eta_H\,\dd r\dd z\right)
 \|L\|_{L^2(\eta_H\otimes\eta_H)}.
\end{equation}
Here the first factor arises from
\[
 \left(\iint_{\Pp^2}r\bar r\,
 \eta_H(X)\eta_H(\bar X)\,\dd X\dd\bar X\right)^{1/2}
 =\iint_{\Pp}r\eta_H(r,z)\,\dd r\dd z.
\]
By \eqref{eq:rz-bound},
\begin{equation}\label{eq:cap-first-moment}
 \iint_{\Pp}r\eta_H\,\dd r\dd z
 \leq\frac{C}{R_T}\iint_{\Pp}rz\eta_H\,\dd r\dd z
 \leq\frac{CM_\Phi(0)}{R_T}.
\end{equation}
Lemma~\ref{lem:log-Lp}, \eqref{eq:P-upper}, and
\(t\in[T,2T]\) imply
\begin{equation}\label{eq:cap-log-bound}
 \|L\|_{L^2(\eta_H\otimes\eta_H)}\leq C\log T.
\end{equation}
Since \(R_T=\kappa T(\log T)^{-2}\),
\begin{equation}\label{eq:cap-energy-small}
 \Ecal(\eta_H)\leq
 C\frac{(\log T)^3}{\kappa T}
 =o_\kappa(1)\qquad(T\to\infty).
\end{equation}
By \eqref{eq:outer-vorticity-energy},
\eqref{eq:cap-energy-small}, and the triangle inequality for the
corresponding velocities, for all sufficiently large \(T\),
\[
 \sqrt{\Ecal(\eta_{\rm rad})}
 \geq\sqrt{\Ecal(\eta_{\rm out})}-\sqrt{\Ecal(\eta_H)}
 \geq\frac12\sqrt{cE_0}.
\]
Thus \eqref{eq:critical-outer-energy} holds, for instance, with
\(E_*=cE_0/4\).
\end{proof}

With the conclusions established above, we are now in a position to prove Theorem~\ref{thm:escape-collision}.

\begin{proof}[Proof of Theorem~\ref{thm:escape-collision}]
Write \(\eta=-\omega\) on \(\Pp\), and set
\[
 E:=\Ecal(\eta(\cdot,t))
 =\frac12\|u_0\|_{L^2(\R^3)}^2=2\pi E_0.
\]
Fix \(A>0\) and \(0<\delta<1\).  Choose
\(\varepsilon_0>0\) so small that
\begin{equation}\label{eq:collision-epsilon-choice}
 1-3\varepsilon_0>\sqrt{1-\delta}.
\end{equation}
We first repeat the exterior localization above while retaining the
size of the error.  Choose \(0<\vartheta<1\) small and \(C_1>1\)
large enough that
\begin{equation}\label{eq:collision-parameter-choice}
 \sqrt{1-\vartheta}-C_1^{-3/2}>1-\varepsilon_0.
\end{equation}
For this proof, put
\[
 R_T:=2A T(\log T)^{-2}.
\]
Proposition~\ref{prop:critical-inner-energy}, applied with radius
\(C_1R_T\), shows that
\[
 \iint_{\Pp\cap\{\varrho\leq C_1R_T\}}r|u|^2\,\dd r\dd z
 \leq CC_1R_T(\log T)^2\mathfrak D(t)
\]
for all sufficiently large \(T\).  Remove from \([T,2T]\) the times
at which
\[
 \mathfrak D(t)>
 \frac{\vartheta E_0}
 {CC_1R_T(\log T)^2}.
\]
We also remove the common null set on which any of the estimates used
below may fail.  If \(G_T\) denotes the remaining set, then
\begin{align}
 |[T,2T]\setminus G_T|
 \leq
 \frac{CC_1R_T(\log T)^2}{\vartheta E_0}
 \int_T^{2T}\mathfrak D(t)\,\dd t
 =o_{A,\delta}(T).
 \label{eq:collision-good-measure}
\end{align}
 Here \(C_1\) and \(\vartheta\) have already been fixed in terms of
\(\delta\), and we used \(\mathfrak D\in L^1(0,\infty)\).

For \(t\in G_T\), angular integration and symmetry in \(z\) give
\[
 \|u(t)\|_{L^2(B_{C_1R_T})}^2
 \leq\vartheta\|u_0\|_{L^2(\R^3)}^2.
\]
Consequently,
\begin{equation}\label{eq:collision-full-outer}
 \|u(t)\|_{L^2(\R^3\setminus B_{C_1R_T})}
 \geq\sqrt{1-\vartheta}\,\|u_0\|_{L^2(\R^3)}.
\end{equation}
Split
\[
 \eta=\eta_{\rm in}+\eta_{\rm out},\qquad
 \eta_{\rm in}:=\eta\ind_{\{\varrho\leq R_T\}},\qquad
 \eta_{\rm out}:=\eta\ind_{\{\varrho>R_T\}}.
\]
Let \(u_{\rm in}=u_{\eta_{\rm in}}\) and
\(u_{\rm out}=u_{\eta_{\rm out}}\).
The kernel in \eqref{eq:energy-static-exact} is nonnegative.  Thus
\(\Ecal(\eta_{\rm in})\leq\Ecal(\eta)=E\), and in particular
\[
 \|u_{\rm in}\|_{L^2(\R^3)}
 \leq\|u_0\|_{L^2(\R^3)}.
\]
Lemma~\ref{lem:harmonic-tail} and
\eqref{eq:collision-full-outer} now imply
\begin{align}
 \|u_{\rm out}\|_{L^2(\R^3)}
 &\geq
 \|u(t)\|_{L^2(\R^3\setminus B_{C_1R_T})}
 -\|u_{\rm in}\|_{L^2(\R^3\setminus B_{C_1R_T})}\notag\\
 &\geq
 \bigl(\sqrt{1-\vartheta}-C_1^{-3/2}\bigr)
 \|u_0\|_{L^2(\R^3)}\notag\\
 &>(1-\varepsilon_0)\|u_0\|_{L^2(\R^3)}.
 \label{eq:collision-outer-vorticity}
\end{align}

We next remove the part that escaped mainly in the axial direction.
Recall from \eqref{eq:high-cap} that, at the present time,
\[
 \eta_H=\eta\ind_{\{\varrho>R_T,\ r\leq R_T/2\}},\qquad
 \eta_{\rm rad}
 =\eta\ind_{\{\varrho>R_T,\ r>R_T/2\}}.
\]
The estimate
\eqref{eq:cap-energy-small}, now with \(\kappa=2A\), gives
\[
 \sup_{t\in[T,2T]}\Ecal(\eta_H)=o_A(1).
\]
Hence, after increasing \(T\), the triangle inequality for the
generated velocities and \eqref{eq:collision-outer-vorticity} yield
\begin{equation}\label{eq:collision-radial-vorticity}
 \|u_{\eta_{\rm rad}}\|_{L^2(\R^3)}
 >(1-2\varepsilon_0)\|u_0\|_{L^2(\R^3)}
 \qquad(t\in G_T).
\end{equation}

It remains to show that the energy-generating radial part lies close
to \(z=0\), as depicted schematically in
Figure~\ref{fig:energy-localization}(b).  Let
\[
 R:=\frac{R_T}{2}=\frac{AT}{(\log T)^2},\qquad
 s_T:=\frac{\log T}{2(\log T-1)}
\]
and, for a constant \(B>1\) to be chosen, set
\[
 H_T:=BR^{-1/2}(\log T)^{1/2}
 =BA^{-1/2}\frac{(\log T)^{3/2}}{\sqrt T}.
\]
Define the high and low parts by
\[
 \eta_{\rm high}:=\eta_{\rm rad}\ind_{\{z>H_T\}},\qquad
 \eta_{\rm low}:=\eta_{\rm rad}\ind_{\{0<z\leq H_T\}}.
\]
The estimates \eqref{eq:energy-e} and \eqref{eq:e-crude},
H\"older's inequality with exponents \(\log T\) and \(2s_T\), and
Lemma~\ref{lem:log-Lp} give
\begin{equation}\label{eq:collision-high-energy-first}
 \Ecal(\eta_{\rm high})\leq C\log T
 \left(\iint_{\Pp}r^{s_T}\eta_{\rm high}\,\dd r\dd z\right)^{1/s_T}.
\end{equation}
Indeed, \(2s_T=\log T/(\log T-1)\) is the H\"older conjugate
of \(\log T\), while \eqref{eq:P-upper} makes the logarithmic-kernel
norm \(O(\log T)\), uniformly for \(t\in[T,2T]\).

On the support of \(\eta_{\rm high}\), one has \(r>R\), \(z>H_T\), and
\(s_T<1\).  The first mixed-moment bound
\eqref{eq:rz-bound} therefore gives
\begin{align}
 \iint_{\Pp}r^{s_T}\eta_{\rm high}\,\dd r\dd z
 &\leq\frac1{H_T}
 \iint_{\Pp}zr^{s_T}\eta_{\rm high}\,\dd r\dd z\notag\\
 &\leq\frac{R^{s_T-1}}{H_T}
 \iint_{\Pp}rz\eta\,\dd r\dd z
 \leq\frac{M_\Phi(0)R^{s_T-1}}{H_T}.
 \label{eq:collision-high-moment}
\end{align}
Since
\[
 1-\frac1{2s_T}=\frac1{\log T},
\]
substitution of \eqref{eq:collision-high-moment} and the definition of
\(H_T\) into \eqref{eq:collision-high-energy-first} yields
\begin{equation}\label{eq:collision-high-energy}
 \Ecal(\eta_{\rm high})\leq
 C B^{-1/s_T}[R\log T]^{1/\log T}.
\end{equation}
For fixed \(A\), the last factor tends to \(e\).  We may therefore
choose \(B\), depending on \(\delta\) and the initial data, so large
that, for all sufficiently large \(T\),
\[
 \|u_{\eta_{\rm high}}\|_{L^2(\R^3)}
 <\varepsilon_0\|u_0\|_{L^2(\R^3)}.
\]

Since
\(u_{\eta_{\rm rad}}=u_{\eta_{\rm low}}+u_{\eta_{\rm high}}\),
\eqref{eq:collision-radial-vorticity} gives
\[
 \|u_{\eta_{\rm low}}\|_{L^2(\R^3)}
 >(1-3\varepsilon_0)\|u_0\|_{L^2(\R^3)}
 >\sqrt{1-\delta}\,\|u_0\|_{L^2(\R^3)}.
\]
Thus \(\Ecal(\eta_{\rm low})>(1-\delta)E\).
Moreover,
\[
 0\leq\eta_{\rm low}\leq
 \eta\ind_{\{r>AT(\log T)^{-2},\ 0<z\leq H_T\}}.
\]
The positive kernel in \eqref{eq:energy-static-exact} makes
\(\Ecal\) monotone under such nonnegative truncations.  Hence the
larger truncation has energy at least \((1-\delta)E\).
Taking \(C_{A,\delta}=BA^{-1/2}\) and using
\eqref{eq:collision-good-measure} proves
\eqref{eq:escape-collision}.
\end{proof}

\subsection{From outer vorticity energy to radial stretching}
\label{sec:growth}

Proposition~\ref{prop:critical-good-times} supplies, at almost every
good time, vorticity at large cylindrical radius whose induced
velocity has a fixed amount of energy.  The remaining step is to
convert that information into radial stretching.  The key point is to
separate the near-diagonal
part of the energy kernel, which is logarithmic but small at large
radius, from the far part, where the stretching kernel is coercive.

\begin{lemma}\label{lem:critical-outer-stretching}
Fix \(E_*>0\).  There are constants \(C_{\rm out}>0\), \(c>0\), and
\(T_{\rm str}>1\),
depending only on \(E_*\) and the initial data, such that the following
holds.  Suppose that
\begin{equation}\label{eq:outer-stretching-range}
 T\geq T_{\rm str},\qquad C_{\rm out}\log T\leq R\leq T.
\end{equation}
Then, for almost every \(t\in[T,2T]\), if there exists a measurable
function \(f\) with
\(0\leq f\leq-\omega(\cdot,t)\), supported in \(\{r>R\}\), and
\(\Ecal(f)\geq E_*\), then
\begin{equation}\label{eq:critical-outer-stretching}
 P'(t)\geq cR^{1/2}(\log T)^{-3/2}.
\end{equation}
\end{lemma}

\begin{proof}
Fix such a time \(t\) and choose such an \(f\).
Let
\begin{equation}\label{eq:critical-IJ}
 I_f:=\iint_{\Pp^2}e(X,\bar X)f(X)f(\bar X)\,\dd X\dd\bar X,
 \qquad
 J_f:=\iint_{\Pp^2}K(X,\bar X)f(X)f(\bar X)\,\dd X\dd\bar X.
\end{equation}
By \eqref{eq:energy-e}, there is \(c_1>0\), depending only on the
initial data, such that \(I_f\geq c_1E_*\).
Recall from \eqref{eq:Sigma} that
\[
 \Sigma
 =\{(X,\bar X)\in\Pp^2:
 B_-(X,\bar X)\leq C_+(X,\bar X)/16\},
\]
the near region illustrated in
Figure~\ref{fig:kernel-reflections}(b).  Write
\[
 I_f=I_{\rm near}+I_{\rm far}
\]
for the contributions from \(\Sigma\) and
\(\Sigma^c\), respectively.  On the near region,
\eqref{eq:near-GMT}, Cauchy--Schwarz,
\eqref{eq:log-Lp-quantitative} with \(p=2\), and
\eqref{eq:rz-bound} give
\begin{align}
 I_{\rm near}
 &\leq C\iint_{\Pp^2}(z\bar z)^{1/2}L(X,\bar X)
 f(X)f(\bar X)\,\dd X\dd\bar X\notag\\
 &\leq C
 \left(\iint_{\Pp}zf\,\dd r\dd z\right)
 \|L\|_{L^2(f\otimes f)}\notag\\
 &\leq C\frac{\log T}{R}.
 \label{eq:critical-near-small}
\end{align}
Here we used
\begin{equation*}
 \iint_{\Pp}zf\,\dd r\dd z
 \leq R^{-1}\iint_{\Pp}rzf\,\dd r\dd z
 \leq M_\Phi(0)R^{-1}.
\end{equation*}
Choosing \(C_{\rm out}\) large enough in
\eqref{eq:outer-stretching-range}, we get
\begin{equation}\label{eq:critical-far-positive}
 I_{\rm far}=I_f-I_{\rm near}\geq\frac{c_1E_*}{2}.
\end{equation}

Put
\begin{equation}\label{eq:critical-holder-parameters}
 \varepsilon_T:=\frac1{\log T},\qquad
 a_T:=\frac{1-2\varepsilon_T}{3}.
\end{equation}
For large \(T\), the four positive numbers
\(a_T,2a_T,\varepsilon_T,\varepsilon_T\) sum to one, and
\(1/4<a_T<1/3\).  Apply \eqref{eq:far-GMT} with
\(\nu=2a_T\).  Since \(r,\bar r>R\),
\begin{equation}\label{eq:critical-radial-weight}
 (r\bar r)^{(1-2a_T)/2}(z\bar z)^{a_T}
 \leq R^{1-4a_T}
 [(rz)(\bar r\bar z)]^{a_T}.
\end{equation}
Indeed, \(1-4a_T<0\) and \(r\bar r\geq R^2\).
We now apply H\"older to the four factors associated with the
\((rz)(\bar r\bar z)\)-moment, \(K\), \(L\), and \(1\), using the
exponents
\[
 \frac1{a_T},\qquad \frac1{2a_T},\qquad
 \frac1{\varepsilon_T},\qquad \frac1{\varepsilon_T},
\]
respectively.  Their reciprocals sum to
\(3a_T+2\varepsilon_T=1\).  Thus
\begin{align}
 I_{\rm far}
 \leq{}&CR^{1-4a_T}
 \left(\iint_{\Pp^2}(rz)(\bar r\bar z)
 f(X)f(\bar X)\,\dd X\dd\bar X\right)^{a_T}\notag\\
 &\quad{}\times J_f^{2a_T}
 \|L\|_{L^{1/\varepsilon_T}(f\otimes f)}
 \left(\iint_{\Pp^2}f(X)f(\bar X)\,\dd X\dd\bar X
 \right)^{\varepsilon_T}\notag\\
 \leq{}&CR^{1-4a_T}J_f^{2a_T}\log T.
 \label{eq:critical-far-holder}
\end{align}
In obtaining the last inequality we have used
\begin{equation*}
 \iint_{\Pp^2}(rz)(\bar r\bar z)f(X)f(\bar X)\,\dd X\dd\bar X
 =\left(\iint_{\Pp}rzf\,\dd r\dd z\right)^2
 \leq M_\Phi(0)^2,
\end{equation*}
the mass bound
\(\iint_{\Pp}f(r,z)\,\dd r\dd z\leq m_0\), and
Lemma~\ref{lem:log-Lp} with
\(1/\varepsilon_T=\log T\), together with
\eqref{eq:P-upper}.  More precisely, uniformly for \(t\in[T,2T]\),
\[
 \|L\|_{L^{1/\varepsilon_T}(f\otimes f)}
 \leq C\bigl[\varepsilon_T^{-1}+\log(2+P(t))\bigr]
 \leq C\log T.
\]
The moment factor is at most \(M_\Phi(0)^{2a_T}\), the mass factor is
at most \((m_0^2)^{\varepsilon_T}\).  Since
\(1/4<a_T<1/3\) and \(0<\varepsilon_T<1\), the resulting constant is
independent of \(T\).

Combining \eqref{eq:critical-far-positive} and
\eqref{eq:critical-far-holder} first gives
\[
 J_f^{2a_T}\geq
 cR^{4a_T-1}(\log T)^{-1}.
\]
Taking the \(2a_T\)-th root yields
\[
 J_f\geq
 cR^{(4a_T-1)/(2a_T)}
 (\log T)^{-1/(2a_T)}.
\]
A direct calculation from
\(a_T=(1-2/\log T)/3\) gives
\[
 \frac{4a_T-1}{2a_T}
 =\frac12-\frac{3}{\log T-2},
 \qquad
 \frac1{2a_T}
 =\frac32+\frac{3}{\log T-2}.
\]
Consequently,
\begin{equation}\label{eq:critical-J-exact}
 J_f\geq
 cR^{\,1/2-\frac{3}{\log T-2}}
 (\log T)^{-3/2-\frac{3}{\log T-2}}.
\end{equation}
Because \(R\leq T\),
\[
 R^{-\frac{3}{\log T-2}}\geq
 \exp\!\left(-\frac{3\log T}{\log T-2}\right),
\]
which is at least \(e^{-6}\) once \(\log T\geq4\).  Moreover,
\[
 (\log T)^{-\frac{3}{\log T-2}}
 =\exp\!\left[-\frac{3\log\log T}{\log T-2}\right]\longrightarrow1,
\]
as \(T\to\infty\).  Thus both correction factors are bounded below by
positive absolute constants for large \(T\).
Therefore
\begin{equation}\label{eq:critical-J-lower}
 J_f\geq cR^{1/2}(\log T)^{-3/2}.
\end{equation}
Finally, \(K\geq0\), \(f\leq-\omega(\cdot,t)\), and
\eqref{eq:Pprime-K} give
\[
\begin{aligned}
 P'(t)
 &\geq c\iint_{\Pp^2}
 K(X,\bar X)[-\omega(X,t)][-\omega(\bar X,t)]
 \,\dd X\dd\bar X\\
 &\geq c\iint_{\Pp^2}K(X,\bar X)f(X)f(\bar X)
 \,\dd X\dd\bar X=cJ_f.
\end{aligned}
\]
Together with \eqref{eq:critical-J-lower}, this proves
\eqref{eq:critical-outer-stretching}.
\end{proof}

We now proceed to prove the remaining main theorems.

\begin{proof}[Proof of Theorem~\ref{thm:main}]
Fix \(\kappa>0\).  For \(T\geq T_{\rm loc}(\kappa)\) and \(t\in G_T\),
Proposition~\ref{prop:critical-good-times} provides
\[
 f=[-\omega(\cdot,t)]
 \ind_{\{\varrho>R_T,\ r>R_T/2\}}
\]
with \(\Ecal(f)\geq E_*\), supported in \(\{r>R_T/2\}\).
Since \(\Ecal(f)>0\), \(f\not\equiv0\), and therefore
\(\mathcal R_\omega(t)\geq R_T/2\) for \(t\in G_T\).  Given \(A>0\),
choose \(\kappa=4A\).  Then, for \(t\in[T,2T]\),
\[
 \frac{R_T}{2}=2A\frac{T}{(\log T)^2}
 \geq A\frac{t}{(\log t)^2}.
\]
Thus, up to a null set,
\[
 G_T\subset
 \left\{t\in[T,2T]:
 \mathcal R_\omega(t)\geq A\frac{t}{(\log t)^2}\right\}.
\]
The complement of the set on the right consequently has measure at
most \(|[T,2T]\setminus G_T|=o_\kappa(T)\), which proves
\eqref{eq:typical-support}.

For each fixed \(\kappa\) and all sufficiently large \(T\), this radius
satisfies
\eqref{eq:outer-stretching-range}.  Lemma~\ref{lem:critical-outer-stretching}
applied with \(R=R_T/2\) therefore gives
\begin{equation}\label{eq:critical-Pprime}
 P'(t)\geq c\kappa^{1/2}T^{1/2}(\log T)^{-5/2}
 \qquad\text{for almost every }t\in G_T.
\end{equation}
Since \(P'\geq0\) and \(|G_T|\geq3T/4\),
\begin{equation}\label{eq:three-halves-increment}
\begin{aligned}
 P(2T)-P(T)
 &=\int_T^{2T}P'(t)\,\dd t\\
 &\geq\int_{G_T}P'(t)\,\dd t\\
 &\geq c\kappa^{1/2}
 T^{1/2}(\log T)^{-5/2}|G_T|\\
 &\geq c\kappa^{1/2}
 \frac{T^{3/2}}{(\log T)^{5/2}}.
\end{aligned}
\end{equation}
Here \(c>0\) is independent of \(\kappa\).  Since \(\kappa>0\) is
arbitrary,
\[
 \lim_{T\to\infty}
 \frac{[P(2T)-P(T)](\log T)^{5/2}}{T^{3/2}}=+\infty.
\]
Indeed, given \(M>0\), first choose the fixed \(\kappa\) so that
\(c\kappa^{1/2}\geq M\), and only then let
\(T\geq T_{\rm loc}(\kappa)\) tend to infinity.
Now put \(t=2T\).  Since \(P\) is nonnegative and increasing,
\[
\begin{aligned}
 \frac{P(t)[\log(2+t)]^{5/2}}{(1+t)^{3/2}}
 &\geq
 \frac{[P(2T)-P(T)][\log(2+2T)]^{5/2}}
 {(1+2T)^{3/2}}\\
 &\geq c\,
 \frac{[P(2T)-P(T)](\log T)^{5/2}}{T^{3/2}}
\longrightarrow+\infty,
\end{aligned}
\]
as \(T\to+\infty\).  This proves \eqref{eq:three-halves-main}.

By the definition of \(\mathcal R_\omega(t)\),
we have
\(P(t)\leq m_0\left(\mathcal R_\omega(t)\right)^2\), and
\eqref{eq:three-halves-support} follows.

Next, we prove the linear estimate on large sets of times.  Fix
\(0<\eta<1\), and choose \(0<\varepsilon<1/10\).  Let
\(L_\varepsilon>1\) and let \(\chi\in C^\infty(\R)\) be non-increasing
and satisfy
\[
 0\leq\chi\leq1,\qquad
 \chi(s)=1\quad(s\leq0),\qquad
 \chi(s)=0\quad(s\geq L_\varepsilon),\qquad
 |\chi'|+|\chi''|\leq\varepsilon.
\]
For \(R>0\), put
\[
 \phi_R(r,z):=\frac{z\varrho^2}{R}
 \chi\!\left(\log\frac{\varrho}{R}\right).
\]
Since \(\phi_R=R^{-1}z(r^2+z^2)\) when \(\varrho\leq R\), this is a
smooth compactly supported multiplier.  With \(\chi,\chi'\), and
\(\chi''\) evaluated at \(\log(\varrho/R)\), direct substitution in
\eqref{eq:multiplier} shows that the sum of its four volume integrands
is
\[
 \frac rR(\chi+\chi')u_{\rm rad}^2
 -\frac{\varrho^2}{Rr}
 \left[\left(2+\frac{z^2}{\varrho^2}\right)\chi+\chi'\right]
 u_{\rm tan}^2
 -\frac zR(3\chi'+\chi'')u_{\rm rad}u_{\rm tan}.
\]
Young's inequality bounds this expression from below by
\[
 \frac1R\ind_{\{\varrho\leq R\}}r|u_{\rm rad}|^2
 -C_\varepsilon\frac{\varrho}{r}|u_{\rm tan}|^2
 -\frac{C\varepsilon}{R}
 \ind_{\{R<\varrho<e^{L_\varepsilon}R\}}r|u|^2.
\]
The boundary terms require no sign condition on the cutoff.  Indeed,
\[
 (\phi_R)_r(r,0)=0,\qquad
 (\phi_R)_z(0,z)=\frac{z^2}{R}
 \left[3\chi\!\left(\log\frac zR\right)
 +\chi'\!\left(\log\frac zR\right)\right],
\]
and hence \(|(\phi_R)_z(0,z)|\leq C_\varepsilon z\).  Moreover,
\[
 0\leq\phi_R\leq e^{L_\varepsilon}z\varrho
 \leq2e^{L_\varepsilon}\Phi,
\]
where $\Phi$ is defined in \eqref{eq:Phi}.
It follows from Lemma~\ref{lem:multiplier} and
\eqref{eq:Phi-perfect-square} that, for almost every \(t\),
\[
 \begin{aligned}
 \frac1R\iint_{\{\varrho\leq R\}}r|u_{\rm rad}|^2\,\dd r\dd z
 \leq{}&\mathcal M_{\phi_R}'(t)
 +C_\varepsilon[-M_\Phi'(t)]\\
 &+\frac{C\varepsilon}{R}
 \iint_{\{R<\varrho<e^{L_\varepsilon}R\}}r|u|^2\,\dd r\dd z.
 \end{aligned}
\]
Since
\[
 r\ind_{\{\varrho\leq R\}}\leq R\frac{\varrho}{r},
\]
the tangential part is also controlled by
\eqref{eq:Phi-perfect-square}.  Integrating on
\(I=[t_0,t_1]\), using energy conservation and
\(0\leq\mathcal M_{\phi_R}(t)\leq
2e^{L_\varepsilon}M_\Phi(0)\), gives
\[
 \begin{aligned}
 \int_I\iint_{\{\varrho\leq R\}}r|u|^2\,\dd r\dd z\dd t
 \leq{}&C_\varepsilon R M_\Phi(0)
 +C_\varepsilon R\int_I[-M_\Phi'(t)]\,\dd t\\
 &+C\varepsilon E_0|I|.
 \end{aligned}
\]

Choose \(C_0>1\) so that \(C_0^{-3/2}\leq1/8\).  We first choose
\(\varepsilon\) sufficiently small in terms of \(\eta\), and then
choose \(a_\eta>0\) sufficiently small.  Since \(M_\Phi\) is
nonnegative and non-increasing,
\[
 \int_T^{2T}[-M_\Phi'(t)]\,\dd t
 =M_\Phi(T)-M_\Phi(2T)\longrightarrow0,
\]
as \(T\to\infty\).  The preceding estimate with \(I=[T,2T]\) and
\(R=C_0a_\eta T\) therefore gives, for all sufficiently large \(T\),
\[
 \int_T^{2T}\iint_{\{\varrho\leq C_0a_\eta T\}}
 r|u|^2\,\dd r\dd z\dd t
 \leq\frac{\eta E_0T}{4}.
\]
After deleting the common null set on which the identities used here
may fail, Chebyshev's inequality provides
\(G_{\eta,T}\subset[T,2T]\) such that
\[
 |G_{\eta,T}|\geq(1-\eta)T,\qquad
 \iint_{\{\varrho\leq C_0a_\eta T\}}r|u|^2\,\dd r\dd z
 \leq\frac{E_0}{4}\quad(t\in G_{\eta,T}).
\]

We now reset \(R_T=a_\eta T\).  For \(t\in G_{\eta,T}\), symmetry in \(z\),
angular integration, and \eqref{eq:E0} give
\[
 \|u(t)\|_{L^2(\R^3\setminus B_{C_0R_T})}
 \geq\frac{\sqrt3}{2}\|u(t)\|_{L^2(\R^3)}.
\]
Splitting the vorticity at \(\varrho=R_T\), then combining Lemma~\ref{lem:harmonic-tail},
the positivity of the energy kernel, and the triangle inequality, as
in \eqref{eq:inner-harmonic-tail}--\eqref{eq:outer-vorticity-energy},
we obtain
\[
 \Ecal\!\left([ -\omega(\cdot,t)]
 \ind_{\{\varrho>R_T\}}\right)\geq cE_0.
\]
The estimates \eqref{eq:cap-energy-first}--\eqref{eq:cap-log-bound},
now with this \(R_T\), give
\[
 \Ecal\!\left([ -\omega(\cdot,t)]
 \ind_{\{\varrho>R_T,\ r\leq R_T/2\}}\right)
 \leq C_\eta\frac{\log T}{T}=o_\eta(1).
\]
The triangle inequality for the corresponding velocities thus yields,
after increasing \(T\),
\[
 \Ecal\!\left([ -\omega(\cdot,t)]
 \ind_{\{\varrho>R_T,\ r>R_T/2\}}\right)>0
 \qquad(t\in G_{\eta,T}).
\]
This truncation is therefore nonzero, and hence
\[
 \mathcal R_\omega(t)
 \geq\frac{R_T}{2}
 \geq\frac{a_\eta}{4}t,
 \qquad t\in G_{\eta,T}.
\]
Together with the measure of \(G_{\eta,T}\), this proves the linear
support assertion in Theorem~\ref{thm:main}, with
\(c_\eta=a_\eta/4\).
\end{proof}

\begin{proof}[Proof of Theorem~\ref{thm:general-Lp}]
The full-time estimate follows from two facts: the three-dimensional
volume of the vorticity support is conserved, and its outer radius is
controlled by the radial velocity.  We begin with the volume.
Put
\[
 \mathcal V_0:=|\{x\in\R^3:\boldsymbol\Omega_0(x)\ne0\}|
 =4\pi\iint_{\{(r,z)\in\Pp:\xi_0(r,z)>0\}}
 r\,\dd r\dd z\in(0,\infty).
\]
The transport of \(\xi\) and preservation of \(r\,\dd r\dd z\) show
that the three-dimensional volume of the vorticity support remains
\(\mathcal V_0\).  Hence H\"older's inequality gives, with
\(1/\infty=0\),
\begin{equation}\label{eq:general-volume-transfer}
 \|\boldsymbol\Omega(t)\|_{L^p(\R^3)}
 \geq \mathcal V_0^{1/p-1}
 \|\boldsymbol\Omega(t)\|_{L^1(\R^3)}
 \geq c\|\boldsymbol\Omega(t)\|_{L^1(\R^3)},
 \qquad 1\leq p\leq\infty,
\end{equation}
where one may take \(c=\min\{1,\mathcal V_0^{-1}\}>0\).  In particular,
the constant is uniform in \(p\).

The radial-velocity estimate in \cite[Proposition~A.1]{LimJeong} gives
\[
 \|u^r(t)\|_{L^\infty(\R^3)}
 \leq C\|u(t)\|_{L^2(\R^3)}^{1/3}
 \left\|\frac{\omega(t)}r\right\|_{L^\infty(\R^3)}^{1/2}
 \|r\boldsymbol\Omega(t)\|_{L^1(\R^3)}^{1/6}.
\]
Energy conservation and transport of \(\omega/r\) make the first two
factors time-independent, while odd symmetry and the sign condition give
\[
 \|r\boldsymbol\Omega(t)\|_{L^1(\R^3)}=4\pi P(t),
\]
and therefore
\begin{equation}\label{eq:general-ur}
 \|u^r(t)\|_{L^\infty(\R^3)}\leq CP(t)^{1/6}.
\end{equation}
Applying \eqref{eq:general-ur} along the meridional flow and taking the
essential supremum over trajectories issuing from the vorticity
support, we obtain
\[
 \mathcal R_\omega(t)
 \leq\mathcal R_\omega(0)+C\int_0^tP(s)^{1/6}\,\dd s
 \leq C(1+t)P(t)^{1/6}.
\]
Here we used the monotonicity of \(P\) and \(P(t)\geq P(0)>0\).
Note that
\[
 P(t)\leq
 \mathcal R_\omega(t)\iint_{\Pp}r[-\omega(r,z,t)]\,\dd r\dd z
 =\frac{\mathcal R_\omega(t)}{4\pi}
 \|\boldsymbol\Omega(t)\|_{L^1(\R^3)}.
\]
Then, using
\eqref{eq:general-volume-transfer}, we obtain
\begin{equation}\label{eq:general-P-to-Lp}
 \inf_{1\leq p\leq\infty}
 \|\boldsymbol\Omega(t)\|_{L^p(\R^3)}
 \geq c\frac{P(t)^{5/6}}{1+t}.
\end{equation}
Consequently,
\[
 \frac{[\log(2+t)]^{25/12}}{(1+t)^{1/4}}
 \inf_{1\leq p\leq\infty}
 \|\boldsymbol\Omega(t)\|_{L^p(\R^3)}
 \geq c\left(
 \frac{P(t)[\log(2+t)]^{5/2}}
 {(1+t)^{3/2}}\right)^{5/6}.
\]
Combining this with \eqref{eq:three-halves-main} yields
\eqref{eq:general-Lp-full}.

We now prove \eqref{eq:general-Lp-typical}.  Fix \(\kappa>0\).  For
\(T\geq T_{\rm loc}(\kappa)\) and \(t\in G_T\), take the truncation \(f\)
supplied by
Proposition~\ref{prop:critical-good-times}, as in the proof of
Theorem~\ref{thm:main}.  Thus \(\Ecal(f)\geq E_*\), and \(f\) is
supported in \(\{r>R_T/2\}\), where
\(R_T=\kappa T(\log T)^{-2}\).
Set
\[
 s_T:=\frac{\log T}{2(\log T-1)}.
\]
For all sufficiently large \(T\), one has \(1/2<s_T<1\).
The exponent is chosen just above \(1/2\) so that
\(2s_T=\log T/(\log T-1)\), the H\"older conjugate of the
logarithmic-kernel exponent \(\log T\).
The estimates \eqref{eq:energy-e} and \eqref{eq:e-crude},
H\"older's inequality, and Lemma~\ref{lem:log-Lp} give
\begin{align}
 E_*
 &\leq C\|L\|_{L^{\log T}(f\otimes f)}
 \left(\iint_{\Pp^2}(r\bar r)^{s_T}
 f(X)f(\bar X)\,\dd X\dd\bar X
 \right)^{1/(2s_T)}\notag\\
 &\leq C(\log T)
 \left(\iint_{\Pp}r^{s_T}f(r,z)\,\dd r\dd z\right)^{1/s_T}.
 \label{eq:general-Lp-low-moment-estimate}
\end{align}
Here \eqref{eq:P-upper} and \(t\in[T,2T]\) are used in the last
inequality.  Hence
\begin{equation}\label{eq:general-Lp-low-moment}
 \iint_{\Pp}r^{s_T}f(r,z)\,\dd r\dd z
 \geq c(\log T)^{-s_T}.
\end{equation}
Since \(f\leq-\omega(\cdot,t)\) and \(r>R_T/2\) on its support,
\[
 \frac1{4\pi}\|\boldsymbol\Omega(t)\|_{L^1(\R^3)}
 \geq\iint_{\Pp}r f(r,z)\,\dd r\dd z
 \geq (R_T/2)^{1-s_T}
 \iint_{\Pp}r^{s_T}f(r,z)\,\dd r\dd z.
\]
Since
\[
 s_T-\frac12=\frac1{2(\log T-1)},
\]
\eqref{eq:general-Lp-low-moment} and the preceding estimate imply
\begin{align*}
 \|\boldsymbol\Omega(t)\|_{L^1(\R^3)}
 &\geq c(R_T/2)^{1-s_T}(\log T)^{-s_T}\\
 &=c\kappa^{1/2}
 \frac{T^{1/2}}{(\log T)^{3/2}}
 [(R_T/2)\log T]^{-(s_T-1/2)}.
\end{align*}
For every fixed \(\kappa>0\),
\[
 [(R_T/2)\log T]^{-(s_T-1/2)}\longrightarrow e^{-1/2}.
\]
It follows from \eqref{eq:general-volume-transfer} that, for all
sufficiently large \(T\),
\[
 \inf_{1\leq p\leq\infty}
 \|\boldsymbol\Omega(t)\|_{L^p(\R^3)}
 \geq c\kappa^{1/2}
 \frac{T^{1/2}}{(\log T)^{3/2}},
 \qquad t\in G_T.
\]
For \(t\in[T,2T]\),
\[
 t^{1/2}(\log t)^{-3/2}
 \leq\sqrt2\,T^{1/2}(\log T)^{-3/2}.
\]
Given \(A>0\), first choose a fixed \(\kappa\) so that
\(c\kappa^{1/2}\geq\sqrt2 A\), and then let \(T\to\infty\).
Together with \eqref{eq:critical-good-measure}, this proves
\eqref{eq:general-Lp-typical}.
\end{proof}

\begin{proof}[Proof of Theorem~\ref{thm:patch-Lp}]
We first list, in the form needed below, the patch moment identities
used in \cite{CJ} and \cite[Corollary~1.3]{EY}.
Preservation of \(r\,\dd r\dd z\), odd reflection, and angular
integration give
\begin{equation}\label{eq:patch-identities}
\begin{aligned}
 m_0&=\iint_{D_t}r\,\dd r\dd z,\\
 P(t)&=\iint_{D_t}r^3\,\dd r\dd z,\\
 \|\boldsymbol\Omega(t)\|_{L^p(\R^3)}^p
 &=4\pi\iint_{D_t}r^{p+1}\,\dd r\dd z,
 \qquad 1\leq p<\infty.
\end{aligned}
\end{equation}
In particular,
\begin{equation}\label{eq:patch-L2-P}
 \|\boldsymbol\Omega(t)\|_{L^2(\R^3)}^2=4\pi P(t),
 \qquad
 |\{x\in\R^3:\boldsymbol\Omega(x,t)\ne0\}|
 =4\pi m_0.
\end{equation}
The patch identity also gives
\[
 \|\boldsymbol\Omega(t)\|_{L^\infty(\R^3)}
 =\operatorname*{ess\,sup}_{(r,z)\in D_t}r
 =\mathcal R_\omega(t).
\]
Together with the linear support estimate in Theorem~\ref{thm:main},
this proves the first assertion of Theorem~\ref{thm:patch-Lp}.

We now prove the density-one family
\eqref{eq:patch-Lp-typical}.  Fix \(\kappa>0\).  For
\(T\geq T_{\rm loc}(\kappa)\) and \(t\in G_T\), let
\[
 f=[-\omega(\cdot,t)]
 \ind_{\{\varrho>R_T,\ r>R_T/2\}}.
\]
By Proposition~\ref{prop:critical-good-times},
\(\Ecal(f)\geq E_*\) and \(f\) is supported in
\(\{r>R_T/2\}\), where \(R_T=\kappa T(\log T)^{-2}\).  Put
\[
 s_T:=\frac{\log T}{2(\log T-1)}.
\]
For all sufficiently large \(T\), one has \(1/2<s_T<1\), and the
derivation of \eqref{eq:general-Lp-low-moment}, which uses no patch
structure, gives
\[
 \iint_{\Pp}r^{s_T}f(r,z)\,\dd r\dd z
 \geq c(\log T)^{-s_T}.
\]

Fix \(1\leq p<\infty\).  For large \(T\), \(s_T<p\); since
\(r>R_T/2\) on the support of \(f\),
\[
 \iint_{\Pp}r^pf(r,z)\,\dd r\dd z
 \geq (R_T/2)^{p-s_T}
 \iint_{\Pp}r^{s_T}f(r,z)\,\dd r\dd z.
\]
Because \(f\) is the indicated truncation of the patch,
\(f=r\) on \(\{f>0\}\), and hence
\[
 \|\boldsymbol\Omega(t)\|_{L^p(\R^3)}^p
 \geq4\pi
 \iint_{\Pp}r^pf(r,z)\,\dd r\dd z.
\]
Combining these inequalities yields
\[
 \|\boldsymbol\Omega(t)\|_{L^p(\R^3)}
 \geq c_p(R_T/2)^{1-s_T/p}(\log T)^{-s_T/p}.
\]
Since \(s_T-\frac12=[2(\log T-1)]^{-1}\) and, by the definition above,
\(R_T/2=\frac{\kappa}{2}T(\log T)^{-2}\), the additional factor
\([(R_T/2)\log T]^{-(s_T-1/2)/p}\) satisfies, for every fixed
\(\kappa>0\),
\[
 \lim_{T\to\infty}
 [(R_T/2)\log T]^{-(s_T-1/2)/p}=e^{-1/(2p)}.
\]
Indeed, the preceding lower bound has been factored exactly as
\[
 (R_T/2)^{1-s_T/p}(\log T)^{-s_T/p}
 =
 (R_T/2)^{1-\frac1{2p}}(\log T)^{-\frac1{2p}}
 [(R_T/2)\log T]^{-(s_T-\frac12)/p}.
\]
Since the limit is \(e^{-1/(2p)}\), independently of the fixed
\(\kappa\), the final factor is at least
\(\frac12e^{-1/(2p)}\) for all sufficiently large \(T\), depending on
\(p\) and \(\kappa\).
Therefore
\begin{equation}\label{eq:typical-Lp-on-good-set}
 \|\boldsymbol\Omega(t)\|_{L^p(\R^3)}
 \geq c_p\kappa^{1-\frac1{2p}}
 \frac{T^{1-\frac1{2p}}}
 {(\log T)^{2-\frac1{2p}}},
 \qquad t\in G_T.
\end{equation}
For \(p=\infty\), the identity
\(\|\boldsymbol\Omega(t)\|_{L^\infty(\R^3)}
=\mathcal R_\omega(t)\), together with
\(\mathcal R_\omega(t)\geq R_T/2\), gives the \(p=\infty\) counterpart of
\eqref{eq:typical-Lp-on-good-set}.  Given \(A>0\), choose the fixed
\(\kappa\) sufficiently large in \eqref{eq:typical-Lp-on-good-set}.
Since, for \(t\in[T,2T]\) and \(1\leq p\leq\infty\),
\[
 t^{1-\frac1{2p}}(\log t)^{\frac1{2p}-2}
 \leq
 2T^{1-\frac1{2p}}(\log T)^{\frac1{2p}-2},
\]
\eqref{eq:critical-good-measure} proves
\eqref{eq:patch-Lp-typical}.

H\"older's inequality on the conserved volume of the patch support yields,
uniformly for \(2\leq p\leq\infty\),
\begin{equation}\label{eq:patch-Lp-from-P}
 \|\boldsymbol\Omega(t)\|_{L^p(\R^3)}
 \geq(4\pi)^{1/p}m_0^{1/p-1/2}P(t)^{1/2}.
\end{equation}
Here we used
\[
 \|\boldsymbol\Omega(t)\|_{L^2(\R^3)}
 \leq
 |\operatorname{supp}\boldsymbol\Omega(t)|^{\frac12-\frac1p}
 \|\boldsymbol\Omega(t)\|_{L^p(\R^3)}
\]
and then substituted \eqref{eq:patch-L2-P}.
The prefactor in \eqref{eq:patch-Lp-from-P} is a positive continuous
function of \(1/p\in[0,1/2]\); it therefore has a positive lower bound
depending only on \(m_0\).  This proves the asserted uniformity for
\(2\leq p\leq\infty\).  Together with \eqref{eq:three-halves-main}, it proves
\eqref{eq:patch-Linfty-pointwise} and
\eqref{eq:patch-Lp-pointwise-high}.

For \(1\leq p<2\), we apply \cite[Theorem~1.1]{EY}; in their scalar
sign convention, the initial
vorticity is \(-\omega_0\), which has the same norms as \(\omega_0\).
Its hypotheses are satisfied: boundedness of \(D_0\) and its positive
distance from \(\partial\Pp\) imply that
\(\omega_0\) and \(\omega_0/r\) belong to
\(L^1(\R^3)\cap L^\infty(\R^3)\), while the Biot--Savart velocity
\(u_0\) has finite energy.  The required second radial moment is finite
because
\[
 \iint_\Pi r^2|\omega_0(r,z)|\,\dd r\dd z=2P(0)<\infty.
\]
Hence
\begin{equation}\label{eq:EY-Linfty-upper}
 \|\boldsymbol\Omega(t)\|_{L^\infty(\R^3)}
 \leq C(1+t)^{4/3}.
\end{equation}
Using
\[
 \|\boldsymbol\Omega(t)\|_{L^2(\R^3)}^2
 \leq\|\boldsymbol\Omega(t)\|_{L^p(\R^3)}^p
       \|\boldsymbol\Omega(t)\|_{L^\infty(\R^3)}^{2-p}
\]
together with \eqref{eq:patch-L2-P}, \eqref{eq:three-halves-main}, and
\eqref{eq:EY-Linfty-upper}, we obtain
\[
 \|\boldsymbol\Omega(t)\|_{L^p(\R^3)}
 \geq c_pP(t)^{1/p}
 (1+t)^{-4(2-p)/(3p)}.
\]
Since
\[
 \frac{3}{2p}-\frac{4(2-p)}{3p}
 =\frac{8p-7}{6p},
\]
this is equivalently
\[
 \frac{\|\boldsymbol\Omega(t)\|_{L^p(\R^3)}
 [\log(2+t)]^{5/(2p)}}
 {(1+t)^{(8p-7)/(6p)}}
 \geq c_p
 \left(
 \frac{P(t)[\log(2+t)]^{5/2}}
 {(1+t)^{3/2}}
 \right)^{1/p}.
\]
The right-hand side tends to infinity.  At \(p=1\), the polynomial
exponent is \(1/6\).  This proves
\eqref{eq:patch-Lp-pointwise-low} and completes the proof.
\end{proof}

\begin{remark}\label{rem:scale-comparison}
The logarithmic factors follow from
\eqref{eq:critical-inner-energy}.  Since
\(\int_0^\infty\mathfrak D(t)\,\dd t<\infty\) and the energy in \(B_R\)
is bounded by \(CR(\log R)^2\mathfrak D(t)\), one may take
\(R\sim\kappa T(\log T)^{-2}\) on a set of times of measure
\(T-o_\kappa(T)\), for each fixed \(\kappa>0\).
Estimate \eqref{eq:harmonic-tail} changes \(R\) only by a fixed factor,
while
\eqref{eq:critical-outer-stretching} contributes
\(R^{1/2}(\log T)^{-3/2}\).  Time integration therefore gives
\(\kappa^{1/2}T^{3/2}(\log T)^{-5/2}\).
Here \(\kappa\) is fixed before \(T\to\infty\); its arbitrariness makes
the full-time prefactors diverge and permits arbitrarily large
constants in the density-one bounds, but it does not improve
the logarithmic power. Indeed, choosing \(\kappa=\kappa(T)\to\infty\) would require a decay rate
for \(\int_T^{2T}\mathfrak D(t)\,\dd t\), which is not supplied by the mere integrability \(\mathfrak D\in L^1(0,\infty)\).  Thus the argument accounts for the exponents obtained
above without asserting that they are sharp.
\end{remark}

\bigskip

\subsection*{\large Declarations:}

\subsection*{Acknowledgments}
The authors would like to thank Professors In-Jee Jeong, Quoc-Hung Nguyen, and Yao Yao for helpful discussions. D. Cao and J. Fan were supported by National Key R\&D Program (Grant 2023YFA1010001) and NNSF of China (Grant 12371212).  G. Qin was supported by National Key R\&D Program of China (Grant 2025YFA1018400) and NNSF of China (Grant 12471190).

\subsection*{Author Contribution information}
All authors contributed equally.

\subsection*{Conflict of interest statement} On behalf of all authors, the corresponding author states that there is no conflict of interest.

\subsection*{Data availability statement} All data generated or analysed during this study are included in this published article and its supplementary information files.

\appendix
\section{Cutoff justifications for the mixed moments}
\label{app:mixed-justification}

The algebra producing the two coercive quadratic forms is kept in the
main proof because it is the central new part of the argument.  For
completeness, we record here the regularization at the axis and the
removal of the outer cutoff.

\begin{lemma}\label{lem:Phi-justification}
The formal computation leading to
\eqref{eq:Phi-perfect-square} is valid in the time-integrated sense on
every bounded time interval.
\end{lemma}

\begin{proof}
Recall that \(\varrho=(r^2+z^2)^{1/2}\).  For
\(0<\varepsilon\leq1\), define
\[
 \Phi_\varepsilon(r,z)
 :=\int_0^z\sqrt{r^2+s^2+\varepsilon^2}\,\dd s,\qquad
 \varrho_\varepsilon
 :=\sqrt{r^2+z^2+\varepsilon^2}.
\]
Let \(\chi\in C_c^\infty([0,\infty))\) equal one on \([0,1]\) and
vanish on \([2,\infty)\).  We apply
Lemma~\ref{lem:multiplier} to
\[
 \phi_{\varepsilon,R}(r,z)
 :=\chi(\varrho/R)\Phi_\varepsilon(r,z).
\]
On the region where \(\chi=1\), direct differentiation gives
\[
\begin{aligned}
 (\Phi_\varepsilon)_z&=\varrho_\varepsilon,&
 (\Phi_\varepsilon)_r
 &=r\,\operatorname{arsinh}
   \frac{z}{\sqrt{r^2+\varepsilon^2}},\\
 (\Phi_\varepsilon)_{rz}&=\frac r{\varrho_\varepsilon},\qquad&
 (\Phi_\varepsilon)_{zz}
 -(\Phi_\varepsilon)_{rr}&
 +\frac{(\Phi_\varepsilon)_r}{r}
 =\frac z{\varrho_\varepsilon}
   \left(1+\frac{r^2}{r^2+\varepsilon^2}\right).
\end{aligned}
\]
Consequently, the three bulk coefficients in
\eqref{eq:multiplier}, in the order
\((u^r)^2,(u^z)^2,u^ru^z\), are
\[
 -\frac{z^2+\varepsilon^2}{r\varrho_\varepsilon},
 \qquad
 -\frac r{\varrho_\varepsilon},
 \qquad
 \frac z{\varrho_\varepsilon}
 \left(1+\frac{r^2}{r^2+\varepsilon^2}\right).
\]
The boundary coefficients are
\((\Phi_\varepsilon)_r(r,0)=0\) and
\((\Phi_\varepsilon)_z(0,z)=\sqrt{z^2+\varepsilon^2}\).

We first fix \(\varepsilon\) and let \(R\to\infty\) in the identity
integrated between two fixed times \(t_1<t_2\).  On the cutoff annulus
\(\{R\leq\varrho\leq2R\}\), direct differentiation gives, for fixed
\(\varepsilon\),
\[
 |\Phi_\varepsilon|\leq C R^2,\qquad
 |\nabla\Phi_\varepsilon|\leq C_\varepsilon R\log(2+R),\qquad
 |\nabla^2\Phi_\varepsilon|\leq C_\varepsilon\log(2+R).
\]
Expanding the derivatives of
\(\chi(\varrho/R)\Phi_\varepsilon\), every cutoff coefficient is
\(O_\varepsilon(\log(2+R))\), except for a potentially singular term
of size \(O(R/r)\) multiplying \((u^r)^2\).  The former terms vanish
after integration by \(|u|=O(\varrho^{-2})\), while the latter vanishes
by \(|u^r|=O(r\varrho^{-3})\).  The boundary errors on the axes vanish
by the same far-field decay.  This removes the outer cutoff for each
fixed \(\varepsilon\).  The same estimates show that the tails of the
uncut bulk and boundary terms tend to zero as \(R\to\infty\).

It remains to pass to \(\varepsilon\to0\).  On every fixed ball
and bounded time interval,
\eqref{eq:danchin-axis-control} gives \(|u^r|\leq Cr\), while \(u^z\)
is bounded.  Hence
\[
 \frac{z^2+\varepsilon^2}{r\varrho_\varepsilon}|u^r|^2
 \leq Cr\varrho_\varepsilon,\qquad
 \frac r{\varrho_\varepsilon}|u^z|^2\leq C|u^z|^2,
\]
and the mixed coefficient is bounded by \(2\).  These bounds provide an
integrable majorant on each fixed ball, independent of
\(\varepsilon\).  The passage is also uniform in the far field.  For
\(\varrho\geq2\) and \(0<\varepsilon\leq1\),
\[
 \frac{z^2+\varepsilon^2}{r\varrho_\varepsilon}|u^r|^2
 +\frac r{\varrho_\varepsilon}|u^z|^2+|u^ru^z|
 \leq Cr\varrho^{-5},
\]
whose meridional integral over \(\{\varrho>R_1\}\) is \(O(R_1^{-2})\),
uniformly on bounded time intervals.  Likewise, on the \(z\)-axis,
\[
 \sqrt{z^2+\varepsilon^2}|u^z(0,z,t)|^2\leq Cz^{-3}
 \qquad(z\geq2).
\]
Splitting the bulk and boundary integrals at a fixed \(R_1\), applying
dominated convergence inside, and then sending \(R_1\to\infty\) proves
\eqref{eq:Phi-perfect-square}.  The endpoint moments converge as well,
because the vorticity support remains in a fixed compact set on
\([t_1,t_2]\) and \(\Phi_\varepsilon\to\Phi\) locally uniformly.
\end{proof}

\begin{lemma}\label{lem:Psi-justification}
The formal computation for the weight \(\Psi\) in
\eqref{eq:critical-Psi} is valid, and
it yields the integrated identity
\eqref{eq:critical-integrated-identity}.
\end{lemma}

\begin{proof}
We first check the only delicate point, namely local integrability near
\((r,z)=(0,0)\).  Fix \(1/2<\alpha<1\), and choose
\(p>3/(1-\alpha)\).  The local \(W^{1,p}\)-regularity of the velocity
and Morrey's inequality, together with
\eqref{eq:danchin-axis-control} and \(u^z(r,0,t)=0\), imply on every
bounded time interval
\[
 |u^r(r,z,t)|\leq C_Tr,\qquad
 |u^z(r,z,t)|\leq C_Tz^\alpha.
\]
For the spherical components in
\eqref{eq:spherical-components}, it follows that, when
\(\varrho\leq1\),
\[
 |u_{\rm rad}|\leq C_T\varrho^\alpha,\qquad
 |u_{\rm tan}|\leq C_Tr\varrho^{\alpha-1}.
\]
The expansions in \eqref{eq:Theta-small} give
\[
 w+|d_g|+\frac g\varrho\leq\frac C\varrho
 \qquad(0<\varrho\leq1).
\]
Substitution of the preceding bounds into
\eqref{eq:Qg-spherical} and \eqref{eq:critical-matrix} bounds the worst
bulk integrand by \(C_T\varrho^{2\alpha-1}\).  In meridional polar
coordinates this contributes
\(C_T\varrho^{2\alpha}\,\dd\varrho\dd\theta\), and is therefore locally
integrable.  The axis term is locally integrable as well: it is
\(O(z^{2\alpha})\), since \(u^z(0,z,t)=O(z^\alpha)\) and the coefficient
in \eqref{eq:critical-axis-coefficient} remains bounded as
\(z\to0\).

At infinity,
\[
 zg(\varrho)=O\!\left(
 \frac{\varrho^2}{\operatorname{arsinh}\varrho}\right),\qquad
 \nabla[zg(\varrho)]=O\!\left(
 \frac{\varrho}{\operatorname{arsinh}\varrho}\right),\qquad
 \nabla^2[zg(\varrho)]=O(1).
\]
For fixed \(\varepsilon\), apply Lemma~\ref{lem:multiplier} to
\[
 \chi(\varrho/R)
 [\lambda_0\Phi_\varepsilon+zg(\varrho)],
\]
where \(\Phi_\varepsilon\) and \(\chi\) are as in the proof of
Lemma~\ref{lem:Phi-justification}.  The preceding estimates and
\eqref{eq:far-field-decay} show, exactly as for \(\Phi_\varepsilon\),
that all terms involving derivatives of the outer cutoff vanish as
\(R\to\infty\).  After this limit, local domination near the origin and
the uniform far-field tails established in
Lemma~\ref{lem:Phi-justification} permit
\(\varepsilon\to0\) on every bounded time interval.  The
remaining bulk and boundary terms are exactly those in
\eqref{eq:critical-integrated-identity}; the endpoint moments converge
by the same compact-support argument used above.
\end{proof}

\end{document}